\documentclass{amsart}
\usepackage{amssymb, verbatim, multirow, xcolor, graphicx}
\usepackage{hyperref}
\usepackage{amsmath}

\def\KNP{\mathop{\makebox{$\wedge$}\hspace{-10pt}\bigcirc}}
\def\KNPp{\mathop{\makebox{$\wedge$}\hspace{-9pt}\bigcirc}}
\def\Ric{\mathop{\rm Ric}}
\def\cRic{\mathop{\rm R\i\makebox[0pt]{\raisebox{5pt}{\tiny$\circ$\;\,}}c}}

\def\dist{\mathop{\rm dist}}

\def\Riem{\mathop{\rm Rm}}

\def\Vol{\mathop{\rm Vol}}

\def\be{\begin{eqnarray}}
\def\ee{\end{eqnarray}}
\def\beg{\begin{eqnarray*}}
\def\ees{\end{eqnarray*}}
\def\bel{\begin{aligned}}
\def\eel{\end{aligned}}

\def\XXint#1#2#3{{\setbox0=\hbox{$#1{#2#3}{\int}$ }
\vcenter{\hbox{$#2#3$ }}\kern-.6\wd0}}

\newtheorem{thm}{Theorem}[section]
\newtheorem{lem}[thm]{Lemma}
\newtheorem{cor}[thm]{Corollary}
\newtheorem{prop}[thm]{Proposition}
\theoremstyle{definition}

\newtheorem*{rem}{Remark}

\numberwithin{equation}{section}

\begin{document}

\title{Generalized K\"ahler Taub-NUT Metrics and Two Exceptional Instantons}
\date{August 2017}

\author[B. Weber]{Brian Weber}
\address{Department of Mathematics\\
The University of Pennsylvania \\\newline
209 S. 33rd St.\\
Philadelphia PA, 19123, USA}
\email{brweber@math.upenn.edu}

\maketitle

\begin{abstract}
	We study the one-parameter family of generalized Kahler Taub-NUT metrics (discovered by Donaldson), along with two exceptional Taub-NUT-like instantons, and understand them to the extent that should be sufficient for blow-up and gluing arguments.
	In particular we parameterize their geodesics from the origin, determine curvature fall-off rates and volume growth rates for metric balls, and find blow-down limits.
\end{abstract}

\section{Introduction}

We provide information on the generalized K\"ahler Taub-NUT class of metrics and two related exceptional metrics.
These are complete, scalar-flat metrics, each with two commuting holomorphic Killing fields, on the underlying complex manifold $\mathbb{C}^2$.
Understanding the characteristics of these metrics will be important in numerous contexts, such as blow-up analysis of singularities in the extremal K\"ahler context, and gluing constructions.

The generalized Taub-NUT metrics were discovered by Donaldson \cite{Do2} and further studied by Abreu and Sena-Dias \cite{AS} \cite{S}; the two ``exceptional'' instanton metrics are from \cite{W}.
We examine asymptotics such as curvature fall-off and volume growth, and compute $L^2$ curvature energy in a fairly simple way.
For the non-exceptional cases, we show that curvature fall-off is strictly quadratic, except for the standard Taub-NUT metric where curvature fall-off is cubic (as is well known).
We show that both of the exceptional instantons have infinite $L^2$ energy, and actually have quartic volume growth, despite not being ALE.

We also investigate blowdown limits: in the standard Taub-NUT the blowdown is flat $\mathbb{R}^3$, whereas in the generalized case we never obtain smooth manifolds.
By ``blowdown'' we mean scaling the metric $g_\epsilon=\epsilon^2{}g$, sending $\epsilon\rightarrow0$, and taking a pointed Gromov-Hausdorff limit.
Limits in our case are always unique, but the naive expectation that limits be 3-dimensional is wrong; they may be either 2- or 3-dimensional.
To determine this, we find an explicit expression for the collapsing field near infinity and then determine that on spherical shells the leaf-space of this field also foliates the Hopf tori.
This foliation might be rational or irrational, and the collapsing behavior is related to the fact that $\mathbb{S}^3$ might collapse to either a 2-sphere (perhaps with orbifold points), or to a line segment, depending on whether the collapsing field is rational or irrational.

The two exceptional instantons have even more peculiar blowdowns.
Their volume growth is quartic, but their Gromov-Hausdorff blowdowns are 3-manifolds.
Their limiting metrics both have curvature singularities along the entirety of a 1-dimensional submanifold.

The instantons $(N^4,g_4,J,\mathcal{X}^1,\mathcal{X}^2)$ we consider are scalar-flat toric K\"ahler 4-manifolds with commuting real-holomorphic Killing fields $\mathcal{X}^1$, $\mathcal{X}^2$.
The fields $\mathcal{X}^1$, $\mathcal{X}^2$ give the complex manifold $(N^4,J)=\mathbb{C}\times\mathbb{C}$ one of two symmetry structures: rotation on both factors (the generalized Taub-NUT metrics and exceptional Taub-NUT), and translation on one factor and rotation on the other (the exceptional half-plane metric).
We remark that scalar-flat instantons on $\mathbb{C}\times\mathbb{C}$ with two translational fields are always flat, by Corollary 4.5 of \cite{W}.

Taking the metric quotient by the Killing fields produces a 2-manifold $\Sigma^2$ (with boundary) and a metric $g_\Sigma$; the pair $(\Sigma^2,g_\Sigma)$ is called the {\it metric polytope} associated to the instanton.
All metric and curvature information on $N^4$ is encoded in this polytope.
The Ricci curvature is encoded in the {\it Ricci potentials} and the {\it Ricci pseudo-volume form} defined in Section \ref{SubSectionCurvatureQuantities}.
In the appendix we study $W^-$ and show that it takes a surprisingly simple form on the generalized Taub-NUT metrics:
\begin{eqnarray}
	\begin{aligned}
		&W^+\;=\;0, \quad
		W^-\;=\;-K_\Sigma\left(6|\rho|^{-2}\rho\otimes\rho-2Id_{\bigwedge^-} \right)
	\end{aligned} \label{EqnsWeylTensorComplete}
\end{eqnarray}
where $\rho$ is the manifold's Ricci form.
The computation of $W^+$ is due to Derdzinski \cite{Derd}.
Noteworthy is that $W^-$ has just two distinct eigenvalues instead of three.
On general toric K\"ahler 4-manifolds the Weyl tensor has three distnct eigenvalues; it is {\it only} on the generalized Taub-NUT metrics where it has two.
In the Appendix we prove the following proposition.
\begin{prop}
	Assume $(N^4,g_4,J)$ is any scalar-flat toric K\"ahler 4-manifold, and let $\omega^-=d\varphi^1\wedge{}d\varphi^2+Jd\varphi^1\wedge{}Jd\varphi^2$.
	Then $\omega^-\in\bigwedge^-$ and is an eigenform of the Weyl tensor: $W^-(\omega^-)=2K_\Sigma$.
	
	More particularly if $(N^4,g_4,J)$ is a generalized Taub-NUT, then the Ricci form $\rho\in\bigwedge^-$ is an eigenform of the Weyl tensor: $W^-(\rho)=-4K_\Sigma\rho$, and $W^-$ is given by (\ref{EqnsWeylTensorComplete}).
\end{prop}

Of some interest is a new set of explicit examples of singular metrics we find.
In studying certain generalized blowdown limits in Section \ref{SectionBlowDownLimits}, we find a family of scalar-flat K\"ahler 4-manifolds with metrics that are smooth except at one point, where an irremovable curvature singularity exists.
These are not blowdowns in any usual sense.
The underlying manifolds continue to be smooth $\mathbb{C}^2$.
The metrics (while remaining scalar-flat and K\"ahler) are singular.

\subsection{Description of the K\"ahler Reduction}

Here we sketch out the objects under study.
After outlining the momentum construction, we describe the exceptional half-plane instanton, the generalized Taub-NUTs, and the exceptional Taub-NUT.
A fuller development is in Section \ref{SectionKahlerReduction}.

\subsubsection{The moment description and classification}

From a construction originating in classical mechanics, the infinitesimal symplectomorphisms $\mathcal{X}^1$, $\mathcal{X}^2$ lead to canonical  ``action-angle'' coordinates $(\varphi^1,\varphi^2,\theta_1,\theta_2)$, where the ``angles'' $\theta_1$, $\theta_2$ parametrize the integrated flows of the symmetry fields $\mathcal{X}^1$, $\mathcal{X}^2$ and the ``actions'' $\varphi^1$, $\varphi^2$ parametrize the $\mathcal{X}^1$-$\mathcal{X}^2$ leaves themselves.

Any metric $g_4$ on such an instanton can be written explicitly though unenlighteningly in these coordinates.
Projecting to the $(\varphi^1,\varphi^2)$-plane produces exactly the Riemannian projection onto the metric polytope; this is known as the {\it moment map} $(\varphi^1,\varphi^2):N^4\rightarrow\Sigma^2$.
The image is a polygon in the $\varphi^1$-$\varphi^2$ plane with an inherited metric---in the subjects of our study the ``polygon'' is either the quarter-plane or the half-plane.

In \cite{Do2},\cite{AS},\cite{S},\cite{W}, many instantons were studied, in part, by using a set of coordinates called {\it volumetric normal coordinates}.
Denoting these coordinates $(x,y)$, we first set
\begin{eqnarray}
	x=\sqrt{|\mathcal{X}^1|^2|\mathcal{X}^2|^2-\left<\mathcal{X}^1,\mathcal{X}^2\right>^2},
\end{eqnarray}
which is the parallelogram area of $\{\mathcal{X}^1,\mathcal{X}^2\}$.
Remarkably, when $N^4$ is scalar flat, this function is harmonic in the natural polytope metric $g_\Sigma$.
Then $y$ is defined as the harmonic conjugate of $x$, meaning a solution of $dy=-*dx$.
The map $z=x+\sqrt{-1}y$ into $\mathbb{C}$ is analytic, and if the polytope has connected boundary, it is an unbranched map onto the closed right half-plane $\overline{H^2}\subset\mathbb{C}$.
The polytope boundary $\partial\Sigma^2$ maps bijectively onto the imaginary axis.
Then $\varphi^1$, $\varphi^2$ can be expressed in terms of $x$, $y$, and the metrics $g_\Sigma$ and $g_4$ can be written down explicitly in terms of the transition functions.
Indeed in $(x,y)$-coordinates, the polytope metric is simply
\begin{eqnarray}
	g_{\Sigma{}}
	\;=\;\frac{1}{x}\det\left(
	\begin{array}{cc}
	\frac{\partial\varphi^1}{\partial{x}} & \frac{\partial\varphi^2}{\partial{x}} \\
	\frac{\partial\varphi^1}{\partial{y}} & \frac{\partial\varphi^2}{\partial{y}} \\
	\end{array}\right)\left( dx\otimes{}dx\,+\,dy\otimes{}dy \right).
	\label{EqnPolytopeMetric}
\end{eqnarray}
The moment variables $\varphi^1$, $\varphi^2$ are now functions of $(x,y)$, where they are constrained by the degenerate-elliptic PDE
\begin{eqnarray}
	x\left(\varphi^i_{xx}+\varphi^i_{yy}\right)\,-\,\varphi_x\;=\;0. \label{EqnMomPDE}
\end{eqnarray}
In \cite{W} a Liouville-type theorem was used to classify pairs $(\varphi^1,\varphi^2)$ of solutions to this degenerate-elliptic system, under the condition that the corresponding polytope be closed and have connected boundary.
When the polytope is the quarter-plane, it was found in \cite{W} that the only possible metrics are the generalized Taub-NUT metrics first written down in \cite{Do2} and the exceptional Taub-NUT metric.

In the case that the polytope is the half-plane, the metric must be either the flat metric on $\mathbb{C}\times\mathbb{C}$, or a multiple of the exceptional half-plane instanton.


\subsubsection{The exceptional half-plane instanton}

This is the case the polytope is a half-plane, which we may take to be $\Sigma^2=\{\varphi^1\ge0\}$.
As a complex manifold $N^4$ is $\mathbb{C}\times\mathbb{C}$ and the holomorphic actions $\mathcal{X}^1$, $\mathcal{X}^2$ are rotational and translational, respectively.
It was proven in \cite{W} that, after possible affine recombination, the one-parameter family of solutions
\begin{eqnarray}
	\varphi^1=\frac12x^2, \quad \varphi^2=y\,+\,Myx^2
\end{eqnarray}
where $M\ge0$ are the only solutions of (\ref{EqnMomPDE}) that produce the half-plane polytope.
From (\ref{EqnPolytopeMetric}) the corresponding polytope metric is $g_\Sigma=\frac12\left(1+\frac{M}{2}x^2\right)\left((dx)^2+(dy)^2\right)$.
The parameter $M$ simply scales the metric, as seen by the coordinate change $x\mapsto{}x/\sqrt{M}$, $y\mapsto{}y/\sqrt{M}$.
The resulting 4-dimensional instanton is called the exceptional half-plane instanton.
Its full metric $g_4=g_\Sigma+G^{ij}d\theta_i\otimes{}d\theta_j$ and a description of its properties are given in Section \ref{SectionHalfPlaneInstanton}.
The choice $M=0$ produces the flat metric.

\subsubsection{The generalized and exceptional Taub-NUT metrics}

This is the case of the quarter-plane polytope, which we may take to be the first quadrant: $\Sigma^2=\{\varphi^1\ge0,\,\varphi^2\ge0\}$.
The underlying complex manifold is $N^4=\mathbb{C}\times\mathbb{C}$ with two rotational symmetry fields.
After possible affine recombination of $\varphi^1$, $\varphi^2$ there is a precisely two-parameter family of solutions
\begin{eqnarray}
	\begin{aligned}
		&\varphi^1=\frac{1}{\sqrt{2}}\left(-y+\sqrt{x^2+y^2}\right)\,+\,\frac{\alpha}{2}x^2,
		\quad \alpha\ge0 \\
		&\varphi^2=\frac{1}{\sqrt{2}}\left(y+\sqrt{x^2+y^2}\right)\,+\,\frac{\beta}{2}x^2,
		\quad\;\;\; \beta\ge0
	\end{aligned}
\end{eqnarray}
that produces this polytope \cite{W}.
These were written down by Donaldson \cite{Do2} in slightly different coordinates.
The corresponding toric 4-manifolds are the generalized Taub-NUT instantons.
If we set $M=\frac{\alpha+\beta}{2\sqrt{2}}$ and $k=\frac{\alpha-\beta}{\alpha+\beta}$, the polytope metric is
\begin{eqnarray}
	g_\Sigma
	\;=\;\frac{1+2M\left(k\,y+\sqrt{x^2+y^2}\right)}{\sqrt{x^2+y^2}}\left(dx\otimes{d}x+dy\otimes{d}y\right).
\end{eqnarray}
The parameter $M\ge0$ is just scale, as can be seen by the coordinate change $x\mapsto{}x/M$, $y\mapsto{}y/M$.
The choice $M=0$ gives the flat metric on $\mathbb{C}^2$, and the choice $M=1$ gives the standard scale, where $\sup_{N^4}|\sec|=1$.
The parameter $k\in[-1,1]$ is called the instanton's {\it chirality number}, and parametrizes the family of inequivalent Taub-NUT metrics.

The instantons given by $k$ and $-k$ are isometric and the corresponding polytope metrics on $\Sigma^2$ are enantiometric (isometric but with flipped orientation), as seen by simply exchanging the two momentum coordinates.
The case $k=0$ is the standard Taub-NUT (which is achiral and Ricci-flat), and the extreme case $k=1,-1$ is the exceptional Taub-NUT, whose properties are qualitatively different from the other Taub-NUTs.
Choices of $k$ outside the $[-1,1]$ range produce topological and curvature singularities.

\subsection{Description of Results}

The first step is choosing better isothermal coordinates on $\Sigma^2$.
We change to {\it quadratic normal coordinates} $u$, $v$ via the fourth-degree polynomial transitions
\begin{eqnarray}
	\varphi^1=\frac{v^2}{\sqrt{2}M}\left(1+(1+k)u^2\right), \quad\quad
	\varphi^2=\frac{u^2}{\sqrt{2}M}\left(1+(1-k)v^2\right),
\end{eqnarray}
which is a diffeomorphism of the first quadrant to itself.
The $N^4$ metric expressed in $(u,v,\theta_1,\theta_2)$ coordinates is
\begin{eqnarray}
	\begin{aligned}
		g_4&\;=\;\frac{2}{M}\left(1+(1+k)u^2+(1-k)v^2\right)\,\left((du)^2+(dv)^2\right) \\
		&\quad\quad\quad\quad\quad\quad\quad\quad\quad\quad\quad\quad
		\;+\; G^{ij}d\theta_i\otimes{d}\theta_j 
	\end{aligned}\label{EqnNewCoordsMetric}
\end{eqnarray}
where the matrix $(G^{ij})$ is a function of $u$, $v$ the parameter $k$, and the scale factor $M$, but whose particular form is unimportant just now (we write it down in (\ref{EqnMetricInUV})).
In these coordinates we see $M$ explicitly as a scale parameter.

In Section \ref{SectionBlowDownLimits} where we explore blowdowns of our metrics, we find that the symplectomorphic Killing field $\mathcal{X}=(1-k)\mathcal{X}^1-(1+k)\mathcal{X}^2$ is nearly an eigenvector of $G^{ij}d\theta_i\otimes{d}\theta_j$, and asymptotically is precisely an eigenvector.
The corresponding eigenvalue asymptotically approaches a multiple of $1/M$.
The field $\mathcal{X}$ is the collapsing field at infinity in the sense of Cheeger-Gromov collapsing theory \cite{CG1}.
The three remaining eigenvalues of $g_4$ asymptotically grow linearly with distance, reflecting the fact that the manifold's asymptotic volume growth is cubic.

To find parametrized geodesics from the origin $(u,v)=(0,0)$, we use the form of the metric in (\ref{EqnNewCoordsMetric}) and a separation method to solve the Eikonal equation $|\nabla{S}|=1$ explicitly in a certain variety of cases.
Characteristic curves of any Eikonal equation are geodesics, as any function satisfying $|\nabla{}S|=1$ is a distance function, and we find enough of these characteristics to allow explicit parametrization of all geodesics based at the origin.
From this we explicitly compute the polytope metric in exponential polar coordinates (equation (\ref{EqnMetricInPolar})), and then compute the key asymptotic quantities.

We summarize our results in the following.
The first theorem is well-known and is included for completeness.
\begin{thm}[The standard Taub-NUT \cite{H}]
	These are the metrics of (\ref{EqnNewCoordsMetric}) with $k=0$.
	The collapsing field at infinity is $\mathcal{X}_1-\mathcal{X}_2$.
	These metrics are Ricci-flat and have total curvature
	\begin{eqnarray}
		\int|\Riem|^2\;=\;32\pi^2.
	\end{eqnarray}
	Volume growth of geodesic balls is cubic: $\Vol\,B(R)=O(R^3)$ and curvature decay is cubic: $|\Riem|=O(R^{-3})$.
	Its Gromov-Hausdorff blowdown is flat $\mathbb{R}^3$.
\end{thm}

\begin{rem}
	In some works one sees $\int|\Riem|^2=8\pi^2$ for the Taub-NUT metric.
	The difference is a factor of 4, and is due to a different choice of norms.
	Denoting by $|\Riem|_{op}$ the norm of $\Riem$ as an operator $\Riem:\bigwedge^2\rightarrow\bigwedge^2$ and denoting by $|\Riem|^2_{tensor}$ the standard tensor norm, we have $|\Riem|_{tensor}^2=4|\Riem|_{op}^2$.
	Throughout this paper we choose the tensor norm.
	This issue is discussed again after Proposition \ref{PropCurvNorms} and after Lemma \ref{LemmaRicEigenform}.
\end{rem}

\begin{thm}[The chiral Taub-NUTs] \label{ThmChiTaubNUT}
	These are the metrics (\ref{EqnNewCoordsMetric}) with $k\in(-1,0)\cup(0,1)$.
	The collapsing field at infinity is $(1-k)\mathcal{X}^1-(1+k)\mathcal{X}^2$.
	These manifolds are scalar-flat and half-conformally flat, and have total energy
	\begin{eqnarray}
		\begin{aligned}
			&\int|\Ric|^2\;=\;32\pi^2\frac{k^2}{1-k^2}\,, 
			\quad\quad
			\int|W^-|^2\;=\;32\pi^2\frac{1+k^2}{1-k^2}\,, \\
			&\quad\quad\quad\quad\quad\quad
			\int|\Riem|^2\;=\;32\pi^2\frac{1+3k^2}{1-k^2}.
		\end{aligned}
	\end{eqnarray}
	Volume growth of geodesic balls is precisely cubic: $\Vol\,B(R)=O(R^3)$, and curvature decay is quadratic: $|\Ric|=O(R^{-2})$, $|W^-|=O(R^{-2})$.
	Gromov-Hausdorff blowdowns are non-flat and have a curvature singularity point.
	The limit is either a 3-dimensional stratified orbifold (when $k$ is rational) or is the closed half-plane (when $k$ is irrational).
\end{thm}
\begin{rem}
	More exactly, $\lim_{R\rightarrow\infty}R^{-3}Vol\,B(R)=\frac83\pi^2\frac{1}{\sqrt{2M}}\,\left(\frac{1}{\sqrt{1-k}}+\frac{1}{\sqrt{1+k}}\right)$.
	The computation of the $L^2$ norms of $|\Riem|$ and $|\Ric|$ was done in \cite{S}; we include it for completeness, and also because it follows easily from the computational structure built up in Section \ref{SectionKahlerReduction} and the Appendix.
\end{rem}
\begin{thm}[The exceptional Taub-NUT]
	These are the metrics of (\ref{EqnNewCoordsMetric}) with maximum chirality $k=1$ or $-1$.
	They are smooth, geodesically complete, toric, scalar flat, half-conformally flat, and K\"ahler.
	The $L^2$-norms of both $|\Ric|$ and $|W^-|$ are infinite.
	Growth of geodesic balls is quartic: $\Vol\,B(R)=O(R^4)$.
	The Riemann tensor decays quadratically $|\Ric|,|W^-|=O(R^{-2})$ along all geodesics from the origin, except for the family of geodesic rays that make up a certain totally geodesic codimension 2 holomorphic submanifold containing the origin.
	Along these rays, curvature does not decay: $|\Ric|=O(1)$ and $|W^-|=O(1)$.
	
	The Gromov-Hausdorff blowdown of the exceptional Taub-NUT is a 3-dimensional manifold with a curvature singularity that makes up an unbounded codimension-2 submanifold.
\end{thm}

\begin{thm}[The exceptional half-plane instanton]
	In volumetric normal coordinates this instanton has polytope metric
	\begin{eqnarray}
		g_\Sigma=(1+x^2)\left(dx\otimes{d}x+dy\otimes{d}y\right),\; 0\le{x}<\infty,\;-\infty<y<\infty.
		\label{EqnExcHPMetricClean}
	\end{eqnarray}
	The corresponding instanton is smooth, geodesically complete, toric, scalar flat, half-conformally flat, and K\"ahler.
	The metric has $\int|\Ric|^2=\int|\Riem|^2=\infty$.
	Volume growth is quartic: $\Vol\,B(R)=O(R^4)$.
	The Riemann tensor decays quadratically $|\Ric|,|W^-|=O(R^{-2})$ along all geodesics from the origin, except for the family of geodesic rays that make up a certain totally geodesic codimension 2 holomorphic submanifold containing the origin, along which it has no curvature decay: $|\Ric|=O(1)$ and $|W^-|=O(1)$.
	
	The Gromov-Hausdorff blowdown of the exceptional half-plane instanton is a 3-dimensional manifold with a curvature singularity along an unbounded codimension-2 submanifold.
\end{thm}

\begin{rem}
	On the polytope, the exceptional Taub-NUT metric (\ref{EqnExcHPMetricClean}) is
	\begin{eqnarray}
		g\;=\;(1+x^2)(dx\otimes{dx}+dy\otimes{}dy)
	\end{eqnarray}
	expressed in volumetric normal coordinates, whereas the exceptional half-plane instanton (\ref{EqnNewCoordsMetric}) has metric
	\begin{eqnarray}
		g\;=\;(1+u^2)(du\otimes{}du+dv\otimes{}dv)
	\end{eqnarray}
	expressed in quadratic normal coordinates.
	They have suspiciously similar properties in other ways, such as a complete, totally geodesic, codimension 2 submanifold with no curvature decay.
	One might wonder if they are the same, or perhaps if one is a cover of the other.
	But they are different, as we prove at the end of Section \ref{SectionHalfPlaneInstanton}.
	The following remark indicates there is a relationship between them of a different sort.
\end{rem}

\begin{rem}
	The exceptional Taub-NUT has rays along which curvature does not decay.
	The injectivity radius does not collapse along these rays, so a natural question is whether we can compute the pointed Gromov-Hausdorff limit as some basepoint moves to infinity along such a ray.
	We do this in Section \ref{SubSubSecUnscaledLimit}, and find that this limit is the exceptional half-plane instanton.
\end{rem}

\begin{rem}
	All instantons we consider are simply-connected except the exceptional half-plane instanton, which can take a simply connected form as $\mathbb{C}^2$ or a non-simply connected form as $\mathbb{C}\times\mathbb{R}\times\mathbb{S}^1$.
	Since the Killing field $\mathcal{X}^2$ is translational and in particular is nowhere zero, the form $\mathbb{C}\times\mathbb{R}\times\mathbb{S}^1$ may be obtained by taking a quotient along a discrete translational distance.
\end{rem}

\begin{rem}
	The Taub-NUT metric and its generalizations have long been a source of examples in general relativity \cite{NUT} \cite{Taub} \cite{M} and Riemannian geometry.
	For instance there are the multi-Taub-NUT metrics of Gibbons-Hawking \cite{GH}; these are hyperk\"ahler and in particular Ricci-flat, and many of them are toric and so have moment polytopes.
	The only one with a quarter-plane moment polytope is the standard Taub-NUT.
	
	Page \cite{P} explores Euclidean Taub-NUT metrics with a magnetic anomaly, which are Ricci-flat but not half-conformally flat---in particular they are non-K\"ahler---and have curvature singularities.
	Noriaki-Toshihiro \cite{NT} explore classes of ``generalized Taub-NUT'' and ``extended Taub-NUT'' metrics with torus symmetry.
	Like the examples considered here, some of their examples are half-conformally flat and non-Einstein.
	However none are K\"ahler except the standard one.
\end{rem}

\begin{rem}
	We claim explicit solutions of the geodesic equation, but we should say what is meant by ``explicit.''
	Our separation method lets us write down {\it unparameterized} geodesics with simple algebraic expressions (equation (\ref{EqnUnparametrizedEqn})).
	But the parametrization is given by a non-algebraic expression: one must invert a function of the type $f(x)=x^a+x^b+\log(x)$; see equation (\ref{EqnImplicitF}).
	Near infinity we are able to approximate even the parametrization with a simple algebraic expression to arbitrary closeness; see Section \ref{SubSubSecAsymptoticApprox} and especially Corollary \ref{CorAlmostDistanceEst}.
\end{rem}

\begin{rem}
	In this paper we study toric instantons on $\mathbb{C}\times\mathbb{C}$, but toric scalar-flat metrics on the $\mathcal{O}(-l)$ bundles over $\mathbb{P}^1$ are known as well, and were written down in \cite{AS}.
	There it was shown exactly which of these are Einstein: precisely the multi-Taub-NUT and multi-Eguchi-Hanson metrics.
	Of course in the ALE case the Einstein metrics were already known to Kronheimer \cite{K}.
	It would be interesting to learn more about the K\"ahler non-Einstein metrics on these spaces, and in particular what their asymptotics are.
	A reasonable conjecture is that they are asymptotically identical to the metrics studied in this paper.
\end{rem}

\begin{rem}
	Just as there is an ``exceptional Taub-NUT" metric there should be an ``exceptional Eguchi-Hanson'' metric and the like; one for each total space $O(-l)$.
\end{rem}

\begin{rem}
	We use the phrase ``affine recombination'' several times, but there is a delicate point here, as affine recombination can be done in two similar ways that have an important difference.
	First, given two potential functions $\varphi^1$, $\varphi^2$, we may recombine them, without altering anything important about the manifold, by any constant-coefficient affine transformation
	\begin{eqnarray}
		\left(
		\begin{array}{c} \widetilde\varphi^1 \\ \widetilde\varphi^2 \end{array}
		\right)
		\;=\;
		\left(
		\begin{array}{cc} c_1^1 & c_2^1 \\ c_1^2 & c_2^2 \end{array}
		\right)
		\left(
		\begin{array}{c} \varphi^1 \\ \varphi^2 \end{array}
		\right)
		\;+\;
		\left(
		\begin{array}{c} c^1 \\ c^2 \end{array}
		\right). \label{EqnAffineRescale}
	\end{eqnarray}
	So long as the coefficients are constant and the $2\times2$ matrix is in $GL(2,\mathbb{R})$, we still retain two independent symplectomorphic Killing fields $\widetilde{\mathcal{X}}^i=J\nabla\widetilde\varphi{}^i$.
	The polytope itself has been altered by a translation and a planar $GL(2,\mathbb{R})$ transformation, but nothing about the manifold's or the polytope's metric or curvature has changed, except in its coordinate expression.
	Volumetric normal coordinates $x$, $y$ are then created, and the theory proceeds as usual.
	
	On the other had one might create the isothermal coordinates $x$, $y$ first, and then change the potentials $\varphi^1$, $\varphi^2$ without changing $x$, $y$.
	In this case the metric {\it does} change: it is multiplied by the determinant of the coefficient matrix, as can be seen by equation (\ref{EqnsMetricAndSectionalOnPolytope}).
	The metric is unchanged only if the coefficient matrix is in $SL(2,\mathbb{R})$.

\end{rem}

\section{Overview of K\"ahler reduction} \label{SectionKahlerReduction}

We review the moment construction, which compresses all $(N^4,g_4,J,\mathcal{X}^1,\mathcal{X}^2)$ data into a 2-dimensional metric polytope $(\Sigma^2,g_\Sigma)$.
This section is included in part for the reader's convenience and in part to establish notation; it presents techniques developed in \cite{Ab1} \cite{G} \cite{Do2} \cite{AS} \cite{W} and elsewhere.
Only the discussion of Ricci curvature in Section \ref{SubSectionCurvatureQuantities} is new; in that section we introduce the Ricci potentials and the Ricci pseudo-volume form on the polytope.
Most of this section deals with any K\"ahler reduction $N^4\rightarrow\Sigma^2$ where $N^4$ is scalar-flat; we specialize to the Taub-NUTs in Section \ref{SubSecSpecificMetrics}.

\subsection{Polytope construction} \label{SubSecPolytopeConstr}

By assumption the fields $\mathcal{X}^1$, $\mathcal{X}^2$ are Killing fields, infinitesimal symplectomorphisms, and infinitesimal biholomorphisms.
Because $\mathcal{L}_{\mathcal{X}^i}\omega=0$ where $\omega$ is the K\"ahler form of $(N^4,g_4,J)$, we have functions $\varphi^1$, $\varphi^2$ defined up to additive constant by $d\varphi^i=-i_{\mathcal{X}^i}\omega$, which is the same as $\mathcal{X}^i=J\nabla\varphi^i$.
This provides gradient fields $\nabla\varphi^1$, $\nabla\varphi^2$ that commute, so define integrable leaves which are Lagrangian submanifolds.
Assigning to one leaf a value of $(0,0)$ for $(\theta_1,\theta_2)$, then we can then define $(\theta_1,\theta_2)$ functions on the entire manifold as push-forwards along the $\mathcal{X}^1$, $\mathcal{X}^2$ action.
The construction gives the so-called {\it action-angle coordinates} $(\varphi^1,\varphi^2,\theta_1,\theta_2)$; by construction $\mathcal{X}^i=\frac{\partial}{\partial\theta_i}$.
The moment map is just forgetting the angle coordinates, and gives the {\it moment polytope}
\begin{eqnarray}
	\Sigma^2
	\;\triangleq\;Image\left[\left(\varphi^1,\varphi^2\right):N^4\longrightarrow\mathbb{R}^2\right]
\end{eqnarray}
in the $(\varphi^1,\varphi^2)$ plane.
This map is a submersion except where $\mathcal{X}^1$, $\mathcal{X}^2$ have zeros or are collinear.
If $N^4$ is compact, and also in certain non-compact cases such as the metric considered here, it is well known that the image is a closed polygon, although in our setting this ``polygon'' is a quarter-plane or a half-plane.
Because $\mathcal{X}^1$ and $\mathcal{X}^2$ are also Killing, $\Sigma^2$ inherits a Riemannian metric which is obviously smooth in the interior, and is in fact smooth at the boundary except at corners where it is Lipschitz.
Because $[\nabla\varphi^1,\nabla\varphi^2]=0$, the distribution $\{\nabla\varphi^1,\nabla\varphi^2\}$ in $N^4$ is integrable, and indeed by Lemma \ref{LemmaToricBasics} it is totally geodesic.
The polytope is locally isometrically isomorphic to the completion of any of these leaves.
These leaves are Lagrangian so the $N^4$ complex structure does not pass to $\Sigma^2$. $\Sigma^2$ has its own complex structure, the Hodge star.

\subsection{Metric quantities}

In action-angle coordinates $(\varphi^1,\varphi^2,\theta_1,\theta_2)$ the metrics, complex structures, and symplectic structures on $N^4$ and $\Sigma^2$ are
\begin{eqnarray}
	&&\begin{aligned}
		g_4=\left(\begin{array}{c|c}
		G_{ij} & 0 \\
		\hline
		0 & G^{ij}
		\end{array}\right), \;
		J=\left(\begin{array}{c|c}
		0 & -G^{ij} \\
		\hline
		\;G_{ij} & 0
		\end{array}\right), \;
		\omega_4=\left(\begin{array}{c|c}
		0 & -Id \\
		\hline
		Id & 0
		\end{array}\right),
		\end{aligned} \label{EqnsGJOmegaM}
		 \\
		&&\begin{aligned}
		&g_\Sigma \;=\;G, \quad\quad
		J_\Sigma\;=\;\frac{1}{\sqrt{\mathcal{V}}}
		\left(\begin{array}{cc}
		\left<\mathcal{X}^1,\,\mathcal{X}^2\right> & -|\mathcal{X}^1|^2 \\
		|\mathcal{X}^2|^2 & -\left<\mathcal{X}^1,\,\mathcal{X}^2\right>
		\end{array}\right),
	\end{aligned} \label{EqnJSigma}
\end{eqnarray}
where we have set
\begin{eqnarray}
	\begin{aligned}
		&(G^{ij})\;=\;\left(\left<\mathcal{X}^i,\mathcal{X}^j\right>\right)
		\;=\;\left(\left<\nabla\varphi^i,\nabla\varphi^j\right>\right), \quad
		(G_{ij})\;=\;
		(G^{ij})^{-1}, \\
		&\mathcal{V}\;=\;\det(G^{-1})\;=\;|\nabla\varphi^1|^2|\nabla\varphi^2|^2-\left<\nabla\varphi^1,\nabla\varphi^2\right>^2.
	\end{aligned} \label{EqnsMetricAndDet}
\end{eqnarray}
We remark that $G$ is a Hessian, $G_{ij}=\partial^2\mathcal{U}/\partial\varphi^i\partial\varphi^j$, for a function $\mathcal{U}$ known as the {\it symplectic potential} \cite{G} of $N^4$, although we shall not have occasion to use this fact.

Expressing $g_4$ in holomorphic coordinates allows easy computation of scalar and Ricci curvatures.
The $\varphi^i$ are neither pluriharmonic nor even harmonic on $N^4$ or $\Sigma^2$, but the angle coordinates $\theta_i$ are pluriharmonic on $N^4$ and so can be used to determine holomorphic coordinates $(z_1,z_2)$.
It is possible to compute
\begin{eqnarray}
	\begin{aligned}
		&\frac{\partial}{\partial{z}_i}
		\;=\;\frac12\left(\nabla\varphi^i\,-\,\sqrt{-1}\,\mathcal{X}^i\right) \quad \text{and} \quad
		dz_i\;=\; Jd\theta_i \,+\,\sqrt{-1}\,d\theta_i
	\end{aligned} \label{EqnCxCoords}
\end{eqnarray}
so the Hermitian metric is $h=h^{i\bar\jmath}=\left<\partial/\partial{}z_i,\partial/\partial\bar{z}_j\right>=\frac12G^{ij}$ and $\det\,h^{i\bar\jmath}\;=\;\frac14\mathcal{V}$.

Lastly it is important to express the polytope metric $g_\Sigma$ in volumetric coordinates $(x,y)$.
Letting $A=\left(\frac{\partial\varphi^i}{\partial{x}^j}\right)$ be the coordinate transition matrix, the metric and the polytope sectional curvature in $(x,y)$ coordinates are
\begin{eqnarray}
\begin{aligned}
	&g_\Sigma\;=\;\frac{\det(A)}{x}\,\left(dx\otimes{d}x+dy\otimes{d}y\right), \\
	&K_\Sigma\;=\;-\frac{x}{\det(A)}\left(\frac{\partial^2}{\partial{x}^2}
	+\frac{\partial^2}{\partial{y}^2}\right)\,\log\sqrt{\frac{\det(A)}{x}}.
\end{aligned} \label{EqnsMetricAndSectionalOnPolytope}
\end{eqnarray}
The sectional curvature $K_\Sigma$ is the sectional curvature of the quotient space and also of the $\{\nabla\varphi^1,\nabla\varphi^2\}$ leaves in $N^4$, as the leaves are totally geodesic by Lemma \ref{LemmaToricBasics}.

\subsection{Curvature quantities} \label{SubSectionCurvatureQuantities}

Because $det(h^{i\bar\jmath})=\frac14\mathcal{V}$, the Ricci form and scalar curvature of $(N^4,J,\omega)$ are
\begin{eqnarray}
	\begin{aligned}
		\rho&\;=\;-\sqrt{-1}\partial\bar\partial\log\,\mathcal{V}
		\;=\;\frac12dJd\log\,\mathcal{V}, \\
		s&\;=\;-\triangle\log\,\mathcal{V}.
	\end{aligned} \label{EqnsRicScal}
\end{eqnarray}
The function $s$ is $\mathcal{X}^1$, $\mathcal{X}^2$ invariant and so passes down to $\Sigma^2$, where the $(M^4,g_4)$ equation $s=-\triangle\log\,\mathcal{V}$ becomes the $(\Sigma^2,g_\Sigma)$ equation
\begin{eqnarray}
	\triangle_{\Sigma}\sqrt{\mathcal{V}}\;+\;\frac12s\sqrt{\mathcal{V}}\;=\;0.
\end{eqnarray}
We emphasize that $s$ is not the scalar curvature of $(\Sigma^2,g_\Sigma)$, but the scalar curvature of $(N^4,g_4)$ passed down to $\Sigma^2$.
Consequently when $s=0$ on $N^4$ the function $x=\sqrt{\mathcal{V}}$ is harmonic on $\Sigma$, and it has harmonic conjugate $y$, meaning a solution of $dy=-J_{\Sigma}dx$.
By Section 3 of \cite{W}, if the polytope boundary has one component then the complex variable $z=x+iy$ has no critical points, so it is a global complex coordinate that maps $\Sigma^2$ to the right half-plane.

Next we consider how the Ricci curvature of $N^4$ is encoded in the polytope.
The Lie derivative\footnote{Recall the convention $[\mathcal{D}_1,\mathcal{D}_2]=\mathcal{D}_1\mathcal{D}_2-(-1)^{|\mathcal{D}_1||\mathcal{D}_2|}\mathcal{D}_2\mathcal{D}_1$ for derivations.} is $\mathcal{L}_{\mathcal{X}^i}=[d,i_{\mathcal{X}^i}]=di_{\mathcal{X}^i}+i_{\mathcal{X}^i}d$, and because $J$ and $\log\,\mathcal{V}$ are invariant under the fields $\mathcal{X}_i$, we see from (\ref{EqnsRicScal}) that
\begin{eqnarray}
	i_{\mathcal{X}^i}\rho
	=\frac12\mathcal{L}_{\mathcal{X}^i}\left(Jd\log\,\mathcal{V}\right)
	-\frac12d\left(i_{\mathcal{X}^i}Jd\log\,\mathcal{V}\right)
	=d\left<\nabla\varphi^i,\,\nabla\log\,\mathcal{V}^{\frac12}\right>.
\end{eqnarray}
The two functions $\mathcal{R}^i=\left<\nabla\varphi^i,\,\nabla\log\,\mathcal{}x\right>$ we call the {\it Ricci potentials}.
These are invariant functions so pass down to $\Sigma^2$.
On $N^4$ clearly $\rho=-d\mathcal{R}^1\wedge{d}\theta_1-d\mathcal{R}^2\wedge{d}\theta_2$.
In the scalar-flat case we have $\rho\in\bigwedge^-$, meaning $*\rho=-\rho$, and so
\begin{eqnarray}
	\begin{aligned}
		&|\Ric|^2dVol_4\;=\;-2\rho\wedge\rho\;=\;
		4\,d\mathcal{R}^1\wedge{d}\mathcal{R}^2\wedge{d}\theta_1\wedge{d}\theta_2.
	\end{aligned} \label{EqnRicciVolForm}
\end{eqnarray}
The factor of 2 on the $\rho\wedge\rho$ term is owing to the fact that the tensor norm is twice the usual norm on 2-forms: $|\Ric|^2=2*(\rho\wedge*\rho)$.
The 2-form $d\mathcal{R}^1\wedge{d}\mathcal{R}^2$ makes sense on $\Sigma^2$ and is non-negative; we call it the {\it Ricci pseudo-volume form}.
Unlike the potentials $\mathcal{R}^1$, $\mathcal{R}^2$, the Ricci pseudo-volume form is invariant under affine recombination of coordinates $\varphi^1$, $\varphi^2$.

The final curvature quantity to consider is the Weyl curvature.
Of course $W^+=0$ on any scalar-flat K\"ahler manifold \cite{Derd}.
By Lemma \ref{LemmaRicEigenform} we also know that $|W^-|^2=96K_\Sigma{}^2$.

\subsection{The metrics} \label{SubSecSpecificMetrics}

The generalized Taub-NUT instantons have underlying complex manifold $\mathbb{C}\times\mathbb{C}$ with two rotational symmetry fields $\mathcal{X}^1$, $\mathcal{X}^2$, and after possible affine recombination of $\varphi^1$, $\varphi^2$, the polytope is the first quadrant.
The 2-parameter family of moment functions that generate the quarter-plane polytope is
\begin{eqnarray}
	\begin{aligned}
		&\varphi^1
		=\frac{1}{\sqrt{2}}\left(-y+\sqrt{x^2+y^2}\right)+\frac{\alpha}{2}x^2, \quad \alpha\ge0 \\
		&\varphi^2
		=\frac{1}{\sqrt{2}}\left(y+\sqrt{x^2+y^2}\right)+\frac{\beta}{2}x^2, \quad \beta\ge0.
	\end{aligned} \label{EqnsTwoParamMomentFunctions}
\end{eqnarray}
Using $M=\frac{\alpha+\beta}{2\sqrt{2}}$, $k=\frac{\alpha-\beta}{\alpha+\beta}$, we compute the polytope metric and Gaussian curvature
\begin{eqnarray}
	\begin{aligned}
		&g_\Sigma
		\;=\;\frac{1+2M\left(ky+\sqrt{x^2+y^2}\right)}{\sqrt{x^2+y^2}}
		\left(dx\otimes{}dx\,+\,dy\otimes{}dy\right) \\
		&K_\Sigma
		\;=\;M\,\frac{-1+2Mk\left(y+k\sqrt{x^2+y^2}\right)}{\left(1+2M\left(ky+\sqrt{x^2+y^2}\right)\right)^3}
	\end{aligned} \label{EqnsSecCurcPolytope2}
\end{eqnarray}
using equations (\ref{EqnsMetricAndSectionalOnPolytope}) above.
Changing $M$ simply scales the metric (to see this make the simultaneous change $x\mapsto{}x/M$, $y\mapsto{}y/M$), and choosing $M=0$ gives the flat metric.
The parameter $k\in[-1,1]$, the chirality number, changes the metric structure while leaving, say, $K_\Sigma(0,0)$ equal to $-M$.
Therefore $k$ does no scaling.

The exceptional case is $k=1$ (or equivalently $k=-1$), where we see that the negative $y$-axis retains constant $K_\Sigma(0,-y)=-M$, so there is no curvature fall-off along $\{x=0\}$ (or, when $k=-1$, along $\{y=0\}$).

Finally we consider the exceptional half-plane instanton.
The underlying complex manifold is $\mathbb{C}\times\mathbb{C}$; the holomorphic symmetry field $\mathcal{X}_1$ is rotational and $\mathcal{X}_2$ is translational.
The momentum polytope is the half-plane, and after possible affine recombination of $\varphi^1$, $\varphi^2$ we have
\begin{eqnarray}
	\varphi^1\;=\;\frac12x{}^2, \quad \varphi^2\;=\;y\,+\,{M}yx^2
\end{eqnarray}
for any constant $M\ge0$.
We obtain polytope metric and sectional curvature
\begin{eqnarray}
	\begin{aligned}
		&g_\Sigma
		\;=\;\left(1+Mx{}^2\right)\left(dx\otimes{d}x\,+\,dy\otimes{d}y\right), \\
		&K_\Sigma
		\;=\;M\,\frac{-1+Mx^2}{\left(1+Mx^2\right)^3}.
	\end{aligned}
\end{eqnarray}
Replacing $x,y$ by $\tilde{x}=\sqrt{\frac{M}{2}}x$, $\tilde{y}=\sqrt{\frac{M}{2}}y$ we have $g_\Sigma=\frac{1}{M}\left(1+\tilde{x}^2\right)\left(d\tilde{x}^2+d\tilde{y}^2\right)$ and again we see that $M$ is a scale parameter.

\section{Asymptotics of the generalized Taub-NUT metrics} \label{SectionGeneralizedTN}

Here the momentum polytope $\Sigma^2$ is the closed quarter-plane, and the corresponding instantons are the generalized Taub-NUTs.
In \S\ref{SubSectionCoordChanges} we create the very useful {\it quadratic normal} coordinate system.
In \S\ref{SubsecDistanceFunction} we compute the distance function to the origin, express the metric in geodesic normal coordinates, and write down a usable approximation for the distance function.
In \S\ref{SubSecComputationAsymptotics} we use this data to determine the asymptotics of our manifolds, and in \S\ref{SubSecEnergyComp} we compute the $L^2$ norms of the curvature quantities.

\subsection{Quadratic normal coordinates} \label{SubSectionCoordChanges}

From $x$, $y$ coordinates, we define $u$, $v$ coordinates:
\begin{eqnarray}
	u\;=\;\sqrt{M}\sqrt{\sqrt{x^2+y^2}+y\;}, \quad
	v\;=\;\sqrt{M}\sqrt{\sqrt{x^2+y^2}-y\;}.
\end{eqnarray}
These are indeed isothermal coordinates, and in fact are a complex square root of the $(x,y)$ coordinates: $x+\sqrt{-1}y=-\frac{\sqrt{-1}}{2M}(v+\sqrt{-1}u)^2$.
We call them {\it quadratic normal coordinates}.
The inverse transformation is $x=\frac{1}{M}uv$, $y=\frac{1}{2M}(u^2-v^2)$.

The image of $(u,v)$ from the quarter-plane $\Sigma^2$ is again the quarter-plane, not the right half-plane as it is in $(x,y)$ coordinates.
The moment functions and metric are
\begin{eqnarray}
	\begin{aligned}
		&\varphi^1=\frac{v^2}{\sqrt{2}M}\left(1+(1+k)u^2\right), \quad
		\varphi^2=\frac{u^2}{\sqrt{2}M}\left(1+(1-k)v^2\right), \\
		&g_\Sigma\;=\;\frac{2}{M}\left(1+(1+k)u^2+(1-k)v^2\right)\,\left(du^2\,+\,dv^2\right) \\
		&K_\Sigma\;=\;M\,\frac{-1+k\left((1+k)u^2-(1-k)v^2\right)}{\left(1+(1+k)u^2+(1-k)v^2\right)^3}.
	\end{aligned} \label{EqnsAdoptedMetricSectional}
\end{eqnarray}
Later we shall require use of the full metric in $u,v,\theta_1,\theta_2$ coordinates.
It is
\begin{eqnarray}
	\begin{aligned}
		&\quad\quad\quad\quad
		g_4\;=\;\left(\begin{array}{c|c}
		g_\Sigma & \\
		\hline
		& G^{ij}
		\end{array}\right), \;\; \text{where} \\
		\\
		&(G^{ij})=
		\frac{1}{M}
			\left(\begin{array}{cc}
			\frac{v^2\left(1+2(1+k)u^2+(1+k)^2u^2(u^2+v^2)\right)}{1+(1+k)u^2+(1-k)v^2} & \frac{u^2v^2\left(2+(1-k^2)(u^2+v^2)\right)}{1+(1+k)u^2+(1-k)v^2} \\
			\\
			\frac{u^2v^2\left(2+(1-k^2)(u^2+v^2)\right)}{1+(1+k)u^2+(1-k)v^2} & \frac{u^2\left(1+2(1-k)v^2+(1-k)^2v^2(u^2+v^2)\right)}{1+(1+k)u^2+(1-k)v^2}
		\end{array}\right).
	\end{aligned} \label{EqnMetricInUV}
\end{eqnarray}

\subsection{Distance functions and geodesic normal coordinates} \label{SubsecDistanceFunction}

\subsubsection{The distance functions $S_\eta$}

The form of $g_\Sigma$ in $(u,v)$ coordinates allows a separation of variables technique in finding certain solutions of $|\nabla{S}|=1$.
Supposing $S(u,v)=f(u)+h(v)$ and choosing any parameter $\eta\in[0,\pi/2]$, we use  (\ref{EqnsAdoptedMetricSectional}) to write the equation $|\nabla{}S|^2=1$ as
\begin{eqnarray}
\frac{M}{2}\frac{\left(f_u\right)^2\,+\,\left(h_v\right)^2}{\left(\cos^2\eta+(1+k)u^2\right)+\left(\sin^2\eta+(1-k)v^2\right)}\;=\;1
\end{eqnarray}
which separates into
\begin{eqnarray}
	\begin{aligned}
		&\frac{df}{du}\;=\;\sqrt{\frac{2}{M}}\sqrt{\cos^2\eta+(1+k)u^2}, \\
		&\frac{dh}{dv}\;=\;\sqrt{\frac{2}{M}}\sqrt{\sin^2\eta+(1-k)v^2}.
	\end{aligned}
\end{eqnarray}

\begin{figure}[h!]
	\caption{\it Contour plots of the distance function $S_\eta$ for two values of $\eta$.
	Solid curve is the virtual locus $S_\eta=0$, which touches the polytope only at $(0,0)$.
	Dashed curves are additional level-sets.
	Thin solid curves are characteristics for $S_\eta$, which are geodesics.
	Exactly one characteristic intersects the origin for each $\eta$.
	We have chosen chirality number $k=0.5$.}
	\label{FigSEta}
	\includegraphics[scale=1]{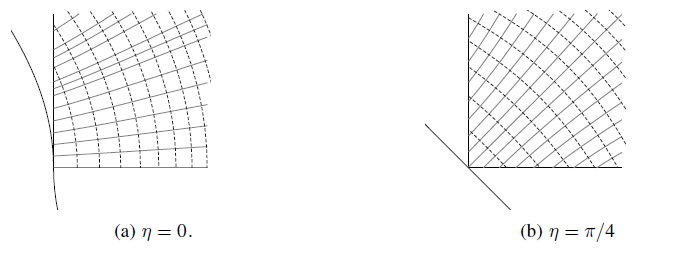}
\end{figure}

Solving for $f$, $h$ under initial conditions $f(0)=h(0)=0$ gives the solution
\begin{eqnarray}
	\begin{aligned}
		S_\eta(u,v)&=
		\sqrt{\frac{2}{M}}\,\frac{\cos^2\eta}{2\sqrt{1+k}}\left[U_\eta\sqrt{1+U_\eta^2}
		\;+\;\log\left(U_\eta+\sqrt{1+U_\eta^2}\,\right)\right] \\
		& \;\;
		+\sqrt{\frac{2}{M}}\,\frac{\sin^2\eta}{2\sqrt{1-k}}\left[V_\eta\sqrt{1+V_\eta^2}
		\;+\;\log\left(V_\eta+\sqrt{1+V_\eta^2}\,\right)\right]
	\end{aligned} \label{EqnDefOfS}
\end{eqnarray}
where we have used the abbreviations $U_\eta=\frac{\sqrt{1+k}}{\cos\eta}u$, $V_\eta=\frac{\sqrt{1-k}}{\sin\eta}v$ and have written $S_\eta$ for $S$ to emphasize the role of the parameter $\eta$.
As depicted in Figure \ref{FigSEta}, the distance function $S_\eta$ is not the distance to any locus within the polytope, but to a virtual locus $S_\eta=0$ in the $u,v$ plane that intersects the polytope only at $(0,0)$.
See Figure \ref{FigSEta}.

\noindent

\subsubsection{The geodesics based at the origin}

Because the virtual locus $S_\eta=0$ intersects the polytope only at the origin, it follows that each choice of $\eta$ allows us to find a single geodesic from the origin.
To study these geodesics from the origin, we attempt to solve for characteristics $\dot\gamma=\nabla{}S_\eta$ with initial condition $\gamma(0,0)=(0,0)$.
For the gradient of $S_\eta$ we have
\begin{eqnarray}
	\begin{aligned}
		\nabla{S}_\eta&\;=\;\sqrt{\frac{2}{M}}\left(
		\frac{\sqrt{\cos^2\eta+(1+k)u^2}}{1+(1+k)u^2+(1-k)v^2}\frac{\partial}{\partial{u}}\right. \\
		&\quad\quad\quad\quad\quad\quad 
		\left.+\frac{\sqrt{\sin^2\eta+(1-k)v^2}}{1+(1+k)u^2+(1-k)v^2}\frac{\partial}{\partial{v}}
		\right).
	\end{aligned}
\end{eqnarray}
With $\gamma(t)=(u(t),v(t))$, the characteristic equation is the coupled autonomous system
\begin{eqnarray}
	\begin{aligned}
		\quad \frac{du}{dt}
		&=\frac{\sqrt{\cos^2\eta+(1+k)u^2}}{1+(1+k)u^2+(1-k)v^2},
		\;\frac{dv}{dt}=\frac{\sqrt{\sin^2\eta+(1+k)v^2}}{1+(1+k)u^2+(1-k)v^2}.
	\end{aligned}
\end{eqnarray}
This is difficult to solve, but eliminating $t$ gives
\begin{eqnarray}
	\frac{dv}{du}\;=\;\frac{\sqrt{\sin^2\eta+(1-k)v^2}}{\sqrt{\cos^2\eta+(1+k)u^2}},
\end{eqnarray}
which separates.
At the point $(u,v)=(0,0)$ we see $\frac{dv}{du}=\tan\eta$, which gives $\eta$ its geometric meaning: it is the initial angle the geodesic makes with the $u$-axis.
The solution for initial condition $\gamma(0)=(0,0)$ is given explicitly by
\begin{eqnarray}
	\left(V_\eta\,+\,\sqrt{1+V_\eta{}^2}\right)^{\frac{1}{\sqrt{1-k}}}
	\;=\;\left(U_\eta\,+\,\sqrt{1+U_\eta{}^2}\right)^{\frac{1}{\sqrt{1+k}}}
	\label{EqnUnparametrizedEqn}
\end{eqnarray}
and again $U_\eta=\frac{\sqrt{1+k}\,u}{\cos\eta}$, $V_\eta=\frac{\sqrt{1-k}\,v}{\sin\eta}$.
This is the {\it unparameterized geodesic equation}.

\subsubsection{Geodesic Normal Coordinates}

\begin{figure}[h!]
	\caption{\it Depictions of geodesic polar coordinates in quadratic normal coordinates and in momentum coordinates.
		Shown are radial geodesics from the origin, and evenly spaced level-sets of the distance function.}
	\label{FigGeodPolar}
	\includegraphics[scale=0.9]{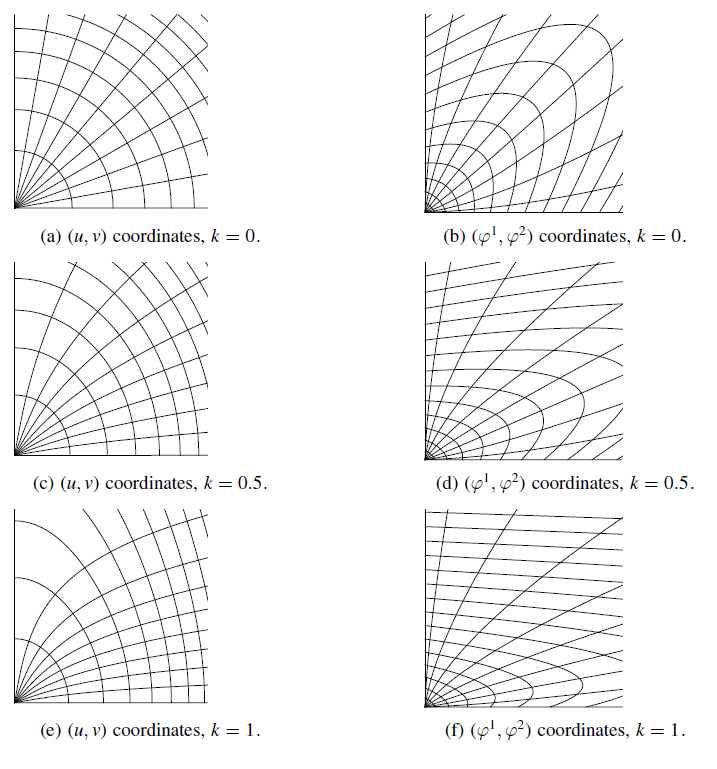}
\end{figure}

Let $R=\dist(o,\cdot)$ be the distance function to the origin.
If $(u,v)$ is an arbitrary point in the first quadrant, we wish to find both the distance $R(u,v)$ to the point and the initial angle $\eta(u,v)$ of the geodesic to that point.
We find the initial angle of the geodesic through a given $(u,v)$ by solving the unparameterized geodesic equation (\ref{EqnUnparametrizedEqn}) for $\eta$:
\begin{eqnarray}
	\begin{aligned}
		&\left(\frac{v\sqrt{1-k}}{\sin\eta}+\sqrt{1+\left(\frac{v\sqrt{1-k}}{\sin\eta}\right)^2}\right)^{\frac{1}{\sqrt{1-k}}} \\
		&\quad\quad\quad\quad
		=\left(\frac{u\sqrt{1+k}}{\cos\eta}+\sqrt{1+\left(\frac{u\sqrt{1+k}}{\cos\eta}\right)^2}\right)^{\frac{1}{\sqrt{1+k}}}.
	\end{aligned} \label{EqnExplicitUnparamGeodEqn}
\end{eqnarray}
This is a non-constructive step.
Given $u$, $v$ there is a unique solution $\eta\in[0,\pi/2]$.
This is because when $\eta$ varies in $[0,\pi/2]$ with $u$ and $v$ fixed, the left-hand side monotonically decreases from $\infty$ and the right-hand side monotonically increases to $\infty$.
Having found $\eta=\eta(u,v)$ this way, the distance to $(u,v)$ is now easy to determine:
\begin{eqnarray}
	R(u,v)\;=\;S_{\eta(u,v)}(u,\,v).
\end{eqnarray}
Thus we have described the transition from the isothermal system $(u,v)$ to polar geodesic coordinates $(R,\eta)$.
This transformation is depicted in Figure \ref{FigGeodPolar}.

We must also compute the reverse transformation: given initial angle $\eta$ and distance $R$, we must find $(u,v)$.
This is equivalent to finding the parametrization for the geodesics described by (\ref{EqnUnparametrizedEqn}).
Given $(R,\eta)$ we must solve the non-algebraic system
\begin{eqnarray}
	\left(U_\eta+\sqrt{1+U_\eta{}^2}\right)^{\frac{1}{\sqrt{1+k}}}=\left(V_\eta+\sqrt{1+V_\eta{}^2}\right)^{\frac{1}{\sqrt{1-k}}}, \; S_\eta(u,v)=R \label{EqnsGeodesicSystem}
\end{eqnarray}
for $(u,v)$.
To do so, we define an auxiliary function $F$ by
\begin{eqnarray}
	F\;\triangleq\;\left(U_\eta+\sqrt{1+U_\eta^2}\right)^{\frac{1}{\sqrt{1+k}}}\;=\;\left(V_\eta+\sqrt{1+V_\eta^2}\right)^{\frac{1}{\sqrt{1-k}}} \label{EqnAuxilliaryFuncDef}
\end{eqnarray}
so that from $F$ we may determine $u$ and $v$:
\begin{eqnarray}
	\begin{aligned}
		&u(F)\;=\;\frac{\cos\eta}{2\sqrt{1+k}}\left(F^{\sqrt{1+k}}-F^{-\sqrt{1+k}}\right), \\
		&v(F)\;=\;\frac{\sin\eta}{2\sqrt{1-k}}\left(F^{\sqrt{1-k}}-F^{-\sqrt{1-k}}\right).
	\end{aligned} \label{EqnsUVByF}
\end{eqnarray}
Using (\ref{EqnDefOfS}) we find that $R=S_\eta(u,v)$ is precisely
\begin{eqnarray}
	\begin{aligned}
		&R
		\;=\;
		\sqrt{\frac{2}{M}}
		\frac{\cos^2\eta}{2\sqrt{1+k}}\left[\frac14\left(F^{2\sqrt{1+k}}-F^{-2\sqrt{1+k}}\right)+\log{F}^{\sqrt{1+k}}\right] \\
		&\quad\quad
		+\sqrt{\frac{2}{M}}\frac{\sin^2\eta}{2\sqrt{1-k}}\left[\frac14\left(F^{2\sqrt{1-k}}-F^{-2\sqrt{1-k}}\right)+\log{F}^{\sqrt{1-k}}\right].
	\end{aligned} \label{EqnImplicitF}
\end{eqnarray}
One then inverts this to obtain $F$ in terms of $R$ and $\eta$, clearly a non-constructive step.
To see that a solution exists and is unique for any choice of $\eta$ and $R$, note that the right-hand side of (\ref{EqnImplicitF}), regarded as a function of $F$, is monotone and has range $(-\infty,\infty)$ as $F$ varies in $(0,\infty)$.
After finding $F=F(R,\eta)$ in this way, (\ref{EqnsUVByF}) gives $u=u(R,\eta)$ and $v=v(R,\eta)$.

In Section \ref{SubSubSecAsymptoticApprox} we find simple, explicit expressions that approximate $u$, $v$, and $F$ as functions of $R$ and $\eta$ with good accuracy.

The coordinates $(R,\eta)$ are, of course, geodesic normal coordinates centered at $(0,0)$.
To compute the metric in this system, consider again the unparameterized geodesic equation (\ref{EqnExplicitUnparamGeodEqn}), which relates $\eta$, $u$, and $v$.
Taking an exterior derivative gives
\begin{eqnarray}
	\frac{du\,+\,u\tan\eta\,d\eta}{\sqrt{\cos^2\eta+(1+k)u^2}}\;=\;
	\frac{dv\,-\,v\cot\eta\,d\eta}{\sqrt{\sin^2\eta+(1-k)v^2}}.
\end{eqnarray}
Since we have $|du|^2=|dv|^2=\frac{M}{2}(1+(1+k)u^2+(1-k)v^2)^{-1}$ and $\left<du,dv\right>=0$, we can isolate $d\eta$ and norm to obtain
\begin{eqnarray}
	&&\quad
	|d\eta|^2\frac{2}{M}\left(u\tan\eta\sqrt{\sin^2\eta+(1-k)v^2}
	+v\cot\eta\sqrt{\cos^2\eta+(1+k)u^2}\right)^2=1.
\end{eqnarray}
Using (\ref{EqnsUVByF}) to write $|d\eta|^2$ in terms of $R$ and $\eta$, we obtain, finally, the polytope metric in geodesic normal coordinates:
\begin{eqnarray}
	\begin{aligned}
		&g_\Sigma\;=\;dR\otimes{d}R\,+\,A(R,\eta)^2d\eta\otimes{d}\eta, \quad \text{where} \\
		&A(R,\eta)^2\;=\;|d\eta|^{-2} \\
		&\quad\quad\quad\;\;=\;\left[\frac{\sin^2\eta}{2M\sqrt{1+k}}\left(F^{\sqrt{1+k}}-F^{-\sqrt{1+k}}\right)\left(F^{\sqrt{1-k}}+F^{-\sqrt{1-k}}\right)\right. \\
		&\quad\quad\quad\quad+\left.\frac{\cos^2\eta}{2M\sqrt{1-k}}\left(F^{\sqrt{1+k}}+F^{-\sqrt{1+k}}\right)\left(F^{\sqrt{1-k}}-F^{-\sqrt{1-k}}\right)\right]^2.
	\end{aligned} \label{EqnMetricInPolar}
\end{eqnarray}

\subsubsection{Asymptotic approximations of $F$ and $R$} \label{SubSubSecAsymptoticApprox}

The coordinates $(u,v)$ and the auxiliary function $F$ are functions of the polar coordinates $(R,\eta)$.
In this section we approximate $u$, $v$ and $F$ using closed-form expressions.
We may approximate the value of $F$ to within tolerable margins by
\begin{eqnarray}
	\begin{aligned}
		&\widetilde{F}(R,\eta)= 
		\begin{cases}
			\left(\frac{8\sqrt{1+k}}{\cos^2\eta}\sqrt{\frac{M}{2}}R\right)^{\frac{1}{2\sqrt{1+k}}},
			& 0\le\eta<\eta_0 \\
			\left(\frac{8\sqrt{1-k}}{\sin^2\eta}\sqrt{\frac{M}{2}}R\right)^{\frac{1}{2\sqrt{1-k}}},
			& \eta_0\;\le\;\eta\;\le\;\frac{\pi}{2},
		\end{cases} \\
		&\eta_0\;=\;\sin^{-1}\left(\frac{\left(\sqrt{\frac{M}{2}}R\right)^{\sqrt{\frac{1+k}{1-k}}-1}}{\left(\sqrt{\frac{M}{2}}R\right)^{\sqrt{\frac{1+k}{1-k}}-1}+\frac{6\sqrt{1+k}}{(8\sqrt{1-k})^{\sqrt{\frac{1+k}{1-k}}}}}\right).
	\end{aligned}
\end{eqnarray}
where ``tolerable margins'' means the following.
\begin{lem} \label{LemmaApproxF}
	Given any values $F$, $\eta$, define the function $\mathcal{R}=\mathcal{R}(F,\eta)$ by
	\begin{eqnarray}
		\begin{aligned}
			&\mathcal{R}\;=\;
			\frac{\cos^2\eta}{\sqrt{2M}\sqrt{1+k}}
			\left[\frac14\left(F^{2\sqrt{1+k}}-F^{-2\sqrt{1+k}}\right)+\log{F}^{\sqrt{1+k}}\right] \\
			&\quad\quad
			+\frac{\sin^2\eta}{\sqrt{2M}\sqrt{1-k}}
			\left[\frac14\left(F^{2\sqrt{1-k}}-F^{-2\sqrt{1-k}}\right)+\log{F}^{\sqrt{1-k}}\right].
		\end{aligned}
	\end{eqnarray}
	If $F=F(R,\eta)$ is the auxiliary function of (\ref{EqnImplicitF}) then of course the distance function is exactly $R=\mathcal{R}(F,\eta)$.
	Given any $\epsilon>0$, then for sufficiently large $R$ we have
	\begin{eqnarray}
		\frac{\mathcal{R}(\widetilde{F}(R,\eta),\,\eta)}{R}
		\;=\;\frac{\mathcal{R}(\widetilde{F}(R,\eta),\,\eta)}{\mathcal{R}(F(R,\eta),\eta)}
		\;\in\;\left[1,\,2+\epsilon\right]
	\end{eqnarray}
	for all $\eta\in[0,\pi/2]$.
\end{lem}
{\it Proof}.
Apply the first derivative test in the parameter $\eta$ to learn that the minimum of $\frac1R\mathcal{R}(\widetilde{F}(R,\eta),\,\eta)$ occurs at the endpoints $\eta=0,\pi/2$, and that the maximum occurs at the discontinuity point, where the left and right limits are different.
Then test these points to learn that the minimum is $1$ and the maximum is a bit bigger than $2$.
\qed

In short, our approximation $\widetilde{F}$ for $F$ gives the correct value of $R$ to within about a factor of 2.
With electronic help, this estimate can be improved with very little use of processing power.
Just a single application of Newton's method---or Householder's method, which is better adapted for this problem---will bring this estimate to within an arbitrarily close multiple of $R$ for large $R$.

This estimate for $F$ is valuable in approximating the metric as given in (\ref{EqnMetricInPolar}), but for determining the key asymptotic ratios a direct approximation of $R$ is preferable.
We create approximations for $(R,\eta)$ with new coordinates $(\widetilde{R},\widetilde\eta)$ that we call {\it almost polar coordinates}, given by
\begin{eqnarray}
	\widetilde{R}
	\;=\;
	\sqrt{\frac{1+k}{2M}}\,u^2
	+\sqrt{\frac{1-k}{2M}}\,v^2,
	\quad\quad
	\widetilde{\eta}
	\;=\;\tan^{-1}\left(\sqrt[4]{\frac{1-k}{1+k}}\frac{v}{u} \right)
	\label{EqnsAlmostPolarQuadTransitions}
\end{eqnarray}
with inverse transitions
\begin{eqnarray}
	u=\sqrt[4]{\frac{2M}{1+k}}\sqrt{\widetilde{R}}\cos\widetilde\eta,
	\quad\quad
	v=\sqrt[4]{\frac{2M}{1-k}}\sqrt{\widetilde{R}}\sin\widetilde\eta. \label{EqnReverseAlmostPolarTrans}
\end{eqnarray}
We shall see that the distance function $R$ and the ``almost distance function'' $\widetilde{R}$ are asymptotically very close together.
Unfortunately the angle $\eta$ and the ``almost angle'' $\widetilde\eta$ are not uniformly close together.
The following estimate is the best we are able to manage for $\widetilde\eta$.
\begin{lem} \label{LemmaAnglesEstimates}
	Assuming $\widetilde{R}$ is sufficiently large compared to $M$ and $k\in(0,1)$, then
	\begin{eqnarray}
		\begin{aligned}
		&\frac{\cos\widetilde\eta}{\cos\eta}
		\;\ge\;\min\left\{\frac{1}{\sqrt{2}},\,
		\left[
		\frac{\left(\sqrt[4]{2M(1-k)}\right)^{\frac{\sqrt{1+k}}{\sqrt{1-k}}}}{\sqrt{1+k}}
		\right]
		\sqrt{\widetilde{R}}^{\frac{\sqrt{1+k}}{\sqrt{1-k}}-1}
		\right\} \\
		&\frac{\sin\widetilde\eta}{\sin\eta}
		\;\ge\;\min\left\{\frac{1}{\sqrt{2}},\,
		\left[
		\frac{\left(\sqrt[4]{2M(1+k)}\right)^{\frac{\sqrt{1-k}}{\sqrt{1+k}}}}{\sqrt{1-k}}
		\right]
		\sqrt{\widetilde{R}}^{\frac{\sqrt{1-k}}{\sqrt{1+k}}-1}
		\right\}.
		\end{aligned} \label{EqnAngleInequalities}		
	\end{eqnarray}
	Indeed it suffices that $\sqrt{M(1-k)}\,\widetilde{R}$ be larger than either of
	\begin{eqnarray}
		\frac12\left(\frac{1}{2^{\sqrt{\frac{1+k}{1-k}}}-1}\right)^{\sqrt{\frac{1-k}{1+k}}}, \quad
		\frac12\left(\frac{1}{2^{\sqrt{\frac{1-k}{1+k}}}-1}\right)^{\sqrt{\frac{1+k}{1-k}}}.
		\label{IneqLargerForTR}
	\end{eqnarray}
\end{lem}
\begin{proof}
	We start with the first inequality.
	If already $\frac{\cos\widetilde\eta}{\cos\eta}\ge\frac{1}{\sqrt{2}}$ then the inequality holds immediately, so assume $\frac{\cos\widetilde\eta}{\cos\eta}<\frac{1}{\sqrt{2}}$.
	Using $\sin^2\widetilde\eta=1-\cos^2\widetilde\eta>1-\frac12\cos^2\eta=\frac12+\frac12\sin^2\eta$ we obtain $\frac{\sin^2\widetilde\eta}{\sin^2\eta}>\frac12\frac{1+\sin^2\eta}{\sin^2\eta}\ge1$.
	Referring to the abbreviations $U_\eta=\frac{\sqrt{1+k}}{\cos\eta}u$ and $V_\eta=\frac{\sqrt{1-k}}{\sin\eta}v$, by using $\frac{\sin^2\widetilde\eta}{\sin^2\eta}>1$ and (\ref{EqnReverseAlmostPolarTrans}) we see
	\begin{eqnarray}
		V_\eta
		\;=\;
		\frac{v\sqrt{1-k}}{\sin\eta}
		\;\ge\;\sqrt[4]{2M(1-k)}\sqrt{\widetilde{R}}.
	\end{eqnarray}
	By (\ref{EqnAuxilliaryFuncDef}) we express $F=\left(V_\eta+\sqrt{1+V_\eta^2}\right)^{\frac{1}{\sqrt{1-k}}}$, and since $V_\eta+\sqrt{1+V_\eta^2}>2V_\eta$ we now have
	\begin{eqnarray}
		F
		\;=\;\left(V_\eta+\sqrt{1+V_\eta^2}\right)^{\frac{1}{\sqrt{1-k}}}
		\;>\;\left(2\sqrt[4]{2M(1-k)}\sqrt{\widetilde{R}}\right)^{\frac{1}{\sqrt{1-k}}}.
	\end{eqnarray}
	As long as we assume $\sqrt[4]{2M(1-k)}\sqrt{\widetilde{R}}$ is larger than either of the expressions in (\ref{IneqLargerForTR}), then we can use this estimate for $F$ to obtain an estimate for $u$.
	Using (\ref{EqnsUVByF}) to express $u$ in terms of $F$, we obtain
	\begin{eqnarray}
		\begin{aligned}
			\frac{u\sqrt{1+k}}{\cos\eta}
			\;=\;F^{\sqrt{1+k}}-F^{-\sqrt{1+k}}
			\;>\;\left(\sqrt[4]{2M(1-k)}\sqrt{\widetilde{R}}\right)^{\frac{\sqrt{1+k}}{\sqrt{1-k}}}
			 \label{EqnFracFirstEst}
		\end{aligned}
	\end{eqnarray}
	Finally using the fact that $u=\widetilde{R}\cos\widetilde\eta$ gives the stated conclusion.
	
	The second inequality proceeds identically, exchanging $v$ for $u$ and so on.
\end{proof}

\begin{cor}[Estimate for the almost distance function] \label{CorAlmostDistanceEst}
	Assuming $\widetilde{R}$ is sufficiently large (as given by (\ref{IneqLargerForTR})), we have
	\begin{eqnarray}
		\widetilde{R}\;<\;
		R
		\;<\;\left(1+\epsilon(\widetilde{R})\right)\widetilde{R}
	\end{eqnarray}
	where $\epsilon(\widetilde{R})\rightarrow0$ as $\widetilde{R}\rightarrow\infty$.
\end{cor}
\begin{proof}
	Substituting the transitions $u=\sqrt[4]{\frac{2M}{1+k}}\sqrt{\tilde{R}}\cos\widetilde\eta$, $v=\sqrt[4]{\frac{2M}{1-k}}\sqrt{\tilde{R}}\sin\widetilde\eta$ in the expression $R=S_{\eta(u,v)}(u,v)$ from (\ref{EqnDefOfS}), we obtain
	\begin{eqnarray*}
		\begin{aligned}
			&R
			\;=\;
			\widetilde{R}\cos^2(\widetilde\eta)
			\sqrt{1+\frac{\cos^2\eta}{(1+k)\widetilde{R}\cos^2\widetilde\eta}}
			+\widetilde{R}\sin^2(\widetilde\eta)
			\sqrt{1+\frac{\sin^2\eta}{(1-k)\widetilde{R}\sin^2\widetilde\eta}} \\
			&\quad\quad +\frac{\cos^2\eta}{\sqrt{2M(1+k)}}
			\log\left[
			\sqrt[4]{2M(1+k)}
			\frac{\sqrt{\widetilde{R}}\cos\widetilde\eta}{\cos\eta}
			\left(1+\sqrt{\frac{\cos^2\eta}{\sqrt{1+k}\widetilde{R}\cos^2\widetilde\eta}+1}\;
			\right)\right] \\
			&\quad\quad +\frac{\sin^2\eta}{\sqrt{2M(1-k)}}
			\log\left[
			\sqrt[4]{2M(1-k)}
			\frac{\sqrt{\widetilde{R}}\sin\widetilde\eta}{\sin\eta}
			\left(1+\sqrt{\frac{\sin^2\eta}{\sqrt{1-k}\widetilde{R}\sin^2\widetilde\eta}+1}\;
			\right)\right].
		\end{aligned}
	\end{eqnarray*}
	By Lemma \ref{LemmaAnglesEstimates} the values $\frac{\sqrt{\widetilde{R}}\cos^2\widetilde\eta}{\cos^2\eta}$ and $\frac{\sqrt{\widetilde{R}}\sin^2\widetilde\eta}{\sin^2\eta}$ both grow like a positive power of $\widetilde{R}$, namely like $\sqrt{\widetilde{R}}{}^{\frac{\sqrt{1+k}}{\sqrt{1-k}}}$ or $\sqrt{\widetilde{R}}{}^{\frac{\sqrt{1-k}}{\sqrt{1+k}}}$, respectively.
	Consequently both logarithms are positive, and so we obtain $\widetilde{R}<R$.
	
	For the upper bound on $R$, using Lemma \ref{LemmaAnglesEstimates} again, we see $\frac{\cos^2\eta}{\sqrt{\widetilde{R}}\cos^2\widetilde\eta}$ and $\frac{\sin^2\eta}{\sqrt{\widetilde{R}}\sin^2\widetilde\eta}$ decay like a power of $\widetilde{R}$; this means the coefficients on $\tilde{R}\cos^2(\widetilde{\eta})$ and $\tilde{R}\cos^2(\widetilde{\eta})$ both approach 1.
	An easy estimate shows the logarithm terms are bounded from above by a definite multiple of $\sqrt{\tilde{R}}\log\widetilde{R}$.
	Thus we conclude
	\begin{eqnarray}
		1\;<\;
		\frac{R}{\;\widetilde{R}\;}
		\;\le\;1+\epsilon(\widetilde{R}).
	\end{eqnarray}
\end{proof}

\subsection{Computation of the asymptotic quantities} \label{SubSecComputationAsymptotics}

We make use of the ``almost polar coordinates'' $(\widetilde{R},\widetilde\eta)$ to compute the key asymptotic ratios of the generalized Taub-NUT instantons.
We note that this section works only for the generalized Taub-NUT metrics because the almost polar coordinates in the two exceptional cases are not given by (\ref{EqnsAlmostPolarQuadTransitions}).
The computations for the two exceptional instantons are deferred to Sections \ref{SectionExceptionalTaubNUT} and \ref{SectionHalfPlaneInstanton}, respectively.

Before computing volumes, we must say a word about the ranges of the coordinates.
Certainly $\widetilde{R}\in[0,\infty)$, $\widetilde\eta\in[0,\pi/2]$.
But the ranges of $\theta_1$, $\theta_2$ are somewhat peculiar:
\begin{eqnarray}
	\theta_1,\theta_2\in[0,\sqrt{8}\pi).
\end{eqnarray}
The ranges for $\theta_1$, $\theta_2$ must be determined through understanding the Delzant gluing construction, where the requirement is that, near the polytope edges, the tori close up to create smooth manifolds without conical singularities.
We examine the situation near a boundary point $(u,0)$ on the $v$-axis.
Consider the 2-manifold determined by fixing $u$ and $\theta_2$, and varying the coordinates $v$ and $\theta_1$.
On this submanifold, the 4-manifold metric (\ref{EqnMetricInUV}) restricts to
\begin{eqnarray}
	\begin{aligned}
		g&\;=\;\left[
		\frac{2}{M}\left(1+(1+k)u^2+(1-k)v^2\right)(dv)^2\right. \\
		&\quad\quad\left.
		+\frac{1}{M}
		\frac{v^2\left(
			\left(1+(1+k)u^2\right)^2+(1+k)^2u^2v^2\right)}{1+(1+k)u^2+(1-k)v^2}\,(d\theta_1)^2
		\right] \\
	\end{aligned}
\end{eqnarray}
for fixed $u$.
Using ``big-$O$'' notation, we write this as
\begin{eqnarray}
	\begin{aligned}
		\quad\quad
		g&\;=\;\frac{2(1+(1+k)u^2)}{M}\left[
		\left(1+O(v^2)\right)(dv)^2
		\,+\,\frac12v^2(1+O(v^2))\,(d\theta_1)^2
		\right] \\
		&\;=\;\frac{2(1+(1+k)u^2)}{M}\left[
		\left(1+O(v^2)\right)(dv)^2
		\,+\,v^2(1+O(v^2))\,\left(d\frac{\theta_1}{\sqrt{2}}\right)^2
		\right].
	\end{aligned}
\end{eqnarray}
Thus, for the central point $(v,\theta_1)=(0,0)$ to be a smooth point rather than a cone point, the variable $\theta_1/\sqrt{2}$ must have range along the circle $[0,2\pi)$, meaning $\theta_1$ ranges along $[0,\sqrt{8}\pi)$.
A similar argument works for the parameterization of $\theta_2$.

The ball $B(S)$ of radius $S$ about the origin is the set of points with radius $R<S$.
Likewise let the {\it almost ball} $AB(S)$ of radius $S$ be
\begin{eqnarray}
	AB(S)\;=\;\left\{\,(\widetilde{R},\widetilde\eta,\theta_1,\theta_2)\in{}N^4\;\;\Big|\;\;
		\widetilde{R}\,<\,S  \;\right\}.
\end{eqnarray}
By Lemma (\ref{CorAlmostDistanceEst}), we have $AB(S)\subset{B}(S)\subset{AB}(S(1+\epsilon))$, where $\lim_{S\rightarrow\infty}\epsilon=0$, and therefore $Vol\,B(S)\;<\;Vol\,AB(S)\;<\;Vol\,B(S(1+\epsilon))$.

\begin{prop}
	If $k\in(-1,1)$, then asymptotic volume growth of balls is cubic:
	\begin{eqnarray}
		\lim_{R\rightarrow\infty}R^{-3}\,\Vol\,B(R)
		\;=\;
		\frac83\pi^2\frac{1}{\sqrt{2M}}\,
		\left(\frac{1}{\sqrt{1-k}}+\frac{1}{\sqrt{1+k}}\right).
	\end{eqnarray}
\end{prop}
\begin{proof}
	In $u,v,\theta_1,\theta_2$ coordinates, we can use (\ref{EqnMetricInUV}) to compute the volume form:
	\begin{eqnarray}
		\begin{aligned}
			&\quad{}dVol
			=\frac{2}{M^2}uv\left(1+(1+k)u^2+(1-k)v^2\right)
			\,du\wedge{d}v\wedge{d}\theta_1\wedge{d}\theta_2.
		\end{aligned}
	\end{eqnarray}
	Transitioning to almost polar coordinates we obtain
	\begin{eqnarray}
		\begin{aligned}
			&dVol\;=\;\frac{2}{M}
			\left(1
			+\sqrt{2M(1+k)}\tilde{R}\cos^2\widetilde\eta
			+\sqrt{2M(1-k)}\tilde{R}\sin^2\widetilde\eta\right) \\
			&\quad\quad\quad\quad\quad\quad\quad\quad\quad\quad
			\cdot\frac{\tilde{R}\cos\widetilde\eta\sin\widetilde\eta}{\sqrt{1-k^2}}
			\,d\tilde{R}\wedge{d}\widetilde\eta\wedge{d}\theta_1\wedge{d}\theta_2.
		\end{aligned} \label{EqnVolInUV}
	\end{eqnarray}
	The ranges for the coordinates are $\widetilde\eta\in[0,\pi/2)$ and $\theta_1,\theta_2\in[0,\sqrt{8}\pi)$.
	Integrating along these ranges and integrating $\tilde{R}$ from $0$ to $S$ gives
	\begin{eqnarray}
		&&\Vol\,AB(S)
		=\frac{\frac83\pi^2S^3}{2M\sqrt{1-k^2}}
		\left[3S^{-1}+\left(\sqrt{2M(1+k)}+\sqrt{2M(1-k)}\right)\right].
	\end{eqnarray}
	Using Lemma \ref{CorAlmostDistanceEst} to approximate balls with almost-balls, we have
	\begin{eqnarray}
		&&\begin{aligned}
			\Vol\,B(S)
			&\le\frac{\frac83\pi^2S^3}{2M\sqrt{1-k^2}}
			\left[3S^{-1}+\left(\sqrt{2M(1+k)}+\sqrt{2M(1-k)}\right)\right] \quad \\
			&\le\Vol\,B(S(1+\epsilon(S))).
		\end{aligned}
	\end{eqnarray}
	so we see that volume growth is indeed cubic when $k\in(-1,1)$.
	Taking the limit,
	\begin{eqnarray}
		\begin{aligned}
			&\lim_{S\rightarrow\infty}S^{-3}\Vol\,B(S)
			\;=\;\frac83\pi^2\;
			\frac{\sqrt{2M(1+k)}+\sqrt{2M(1-k)}}{2M\sqrt{1-k^2}}.
		\end{aligned}
	\end{eqnarray}
\end{proof}

\begin{lem} \label{LemmaKAsymptotics}
	If $K_\Sigma$ is the polytope sectional curvature and $k\ne0,\pm1$, then $K_\Sigma=O(R^{-2})$, except along a single path where $K_\Sigma=O(R^{-3})$.
	In almost polar coordinates,
	\begin{eqnarray}
		\lim_{\widetilde{R}\rightarrow\infty}\widetilde{R}^2{}K_\Sigma
		\;=\;\frac{k}{2}
		\frac{\sqrt{1+k}\cos^2\widetilde\eta-\sqrt{1-k}\sin^2\widetilde\eta}
		{\left(\sqrt{1+k}\cos^2\widetilde\eta+\sqrt{1-k}\sin^2\widetilde\eta\right)^3}.
	\end{eqnarray}
	If $k=0$ then $K_\Sigma=O(R^{-3})$ along all paths to infinity.
\end{lem}
\begin{proof}
	Using the expression (\ref{EqnsAdoptedMetricSectional}) and the transitions to $\widetilde{R}$, $\widetilde\eta$ we obtain
	\begin{eqnarray}
		\begin{aligned}
			K_\Sigma&\;=\;M\,\frac{-1+k\sqrt{2M}\,\widetilde{R}\left(\sqrt{1+k}\cos^2\widetilde\eta
				-\sqrt{1-k}\sin^2\widetilde\eta\right)}{\left(1+\sqrt{2M(1+k)}\,\widetilde{R}\cos^2\widetilde\eta
				+\sqrt{2M(1-k)}\,\widetilde{R}\sin^2\widetilde\eta\right)^3} \\
			&\;=\;
				\frac{k}{2}\frac{1}{{\widetilde{R}}^2}\,\frac{-\frac{1}{k\sqrt{2M}\widetilde{R}}+\sqrt{1+k}\cos^2\widetilde\eta
				-\sqrt{1-k}\sin^2\widetilde\eta}{\left(\frac{1}{\sqrt{2M}\widetilde{R}}+\sqrt{1+k}\,\cos^2\widetilde\eta
				+\sqrt{1-k}\,\sin^2\widetilde\eta\right)^3}.
		\end{aligned}\label{EqnKNewInUV}
	\end{eqnarray}
	Taking a limit, then, we obtain
	\begin{eqnarray}
		\lim_{\widetilde{R}\rightarrow\infty}\widetilde{R}^2K_\Sigma
		\;=\;
		\frac{k}{2}\frac{\sqrt{1+k}\cos^2\widetilde\eta
			-\sqrt{1-k}\sin^2\widetilde\eta}
		{\left(\sqrt{1+k}\cos^2\widetilde\eta
			+\sqrt{1-k}\sin^2\widetilde\eta\right)^3}
	\end{eqnarray}
	Therefore $K_\Sigma=O(\tilde{R}^{-2})=O(R^{-2})$ except along a single path which is the path of constant $\widetilde\eta$ where $\frac{\sin\tilde\eta}{\cos\tilde\eta}=\sqrt[4]{1+k}/\sqrt[4]{1-k}$.
	
	When $k=0$ then (\ref{EqnKNewInUV}) gives $K_\Sigma=O(R^{-3})$ everywhere.
\end{proof}

\begin{lem} \label{LemmaRicAsymptotics}
	Let $(\Sigma^2,g_{\Sigma})$ be a Taub-NUT polytope with $k\in[-1,1]$.
	The Ricci potentials are
	\begin{eqnarray}
		\begin{aligned}
			&\mathcal{R}^1
			\;=\;\frac{1}{\sqrt{2}}\frac{1+(1+k)(u^2+v^2)}{1+(1+k)u^2+(1-k)v^2}, \\
			&\mathcal{R}^2
			\;=\;\frac{1}{\sqrt{2}}\frac{1+(1-k)(u^2+v^2)}{1+(1+k)u^2+(1-k)v^2},
		\end{aligned} \label{EqnSpecificRicciPots}
	\end{eqnarray}
	the norm of Ricci curvature is
	\begin{eqnarray}
		\begin{aligned}
			|\Ric|&\;=\;\frac{4|k|M}{\left(1+(1+k)u^2+(1-k)v^2\right)^2},
		\end{aligned}
	\end{eqnarray}
	and we have
	\begin{eqnarray}
		\quad\quad\quad
		|\Ric|^2dVol_4=\frac{2uv}{M^2}\left(1+(1+k)u^2+(1-k)v^2\right)
		du\wedge{}dv\wedge{}d\theta_1\wedge{}d\theta_2. \label{EqnRic2DVol}
	\end{eqnarray}
\end{lem}
\begin{proof}
	From Section \ref{SubSectionCurvatureQuantities} the Ricci potentials are defined by $\mathcal{R}^i=\left<\nabla\varphi^i,\,\nabla\log\mathcal{}x\right>$; an elementary computation gives (\ref{EqnSpecificRicciPots}).
	Using (\ref{EqnRicciVolForm}), we compute $|\Ric|^2dVol_4$:
	\begin{eqnarray}
		\begin{aligned}
			&|\Ric|^2dVol_4
			\;=\;4\,d\mathcal{R}^1\wedge{}d\mathcal{R}^2\wedge{}d\theta_1\wedge{}d\theta_2 \\
			&\quad
			\;=\;\frac{32k^2uv}{\left(1+(1+k)u^2+(1-k)v^2\right)^3}
			\,du\wedge{}dv\wedge{}d\theta_1\wedge{}d\theta_2.
		\end{aligned} \label{EqnComputedRicSqVol}
	\end{eqnarray}
	Using (\ref{EqnMetricInUV}) we compute $dVol_4$ in $(u,v,\theta_1,\theta_2)$ coordinates:
	\begin{eqnarray}
		dVol_4=\frac{2uv}{M^2}\left(1+(1+k)u^2+(1-k)v^2\right)
		\,du\wedge{}dv\wedge{}d\theta_1\wedge{}d\theta_2
	\end{eqnarray}
	so therefore
	\begin{eqnarray}
		\begin{aligned}
			|\Ric|^2&\;=\;\frac{16k^2M^2}{\left(1+(1+k)u^2+(1-k)v^2\right)^4}.
		\end{aligned}
	\end{eqnarray}
\end{proof}

\begin{prop}[Curvature Decay Rates] \label{PropInstantonCurvDecay}
	In the generic case $k\neq0,\pm1$, we have $|\Ric|,|W^-|=O(R^{-2})$.
	In the case $k=0$ we have $|W^-|=O(R^{-3})$.
\end{prop}
\begin{proof}
	If $W=W^++W^-$ is the Weyl tensor, the computation of $W^-$ from (\ref{EqnWeylTensorNorm}) and the fact that $W^+=0$ gives
	\begin{eqnarray}
		|W|^2\;=\;96|K_\Sigma|^2
	\end{eqnarray}
	so from Lemma \ref{LemmaKAsymptotics} we obtain the claimed $R^{-2}$ decay rate for $|W^-|$.
	From Lemma \ref{LemmaRicAsymptotics}
	\begin{eqnarray}
		|\Ric|\;=\;\frac{4|k|M}{\left(1+(1+k)u^2+(1-k)v^2\right)^2}.
	\end{eqnarray}
	Using Corollary \ref{CorAlmostDistanceEst} and changing to the almost polar coordinates, we obtain
	\begin{eqnarray}
		|\Ric|\;=\;\frac{4kM}{\left(1
		+\sqrt{2M(1+k)}\,\widetilde{R}\cos^2\widetilde\eta
		+\sqrt{2M(1-k)}\,\widetilde{R}\sin^2\widetilde\eta\right)^2}.
	\end{eqnarray}
	When $k$ is not $0,1,-1$ we see that $|\Ric|=O(\tilde{R}^{-2})=O(R^{-2})$.
	We now have that both $|\Ric|$ and $|W^-|$ are $O(R^{-2})$.

	When $k=0$ we have $\Ric=0$ and therefore $K_\Sigma=O(R^{-3})$ gives $|W^-|=O(R^{-3})$.
\end{proof}

\subsection{$L^2$ norms} \label{SubSecEnergyComp}

Using the ``Ricci potentials'' from section \ref{SubSectionCurvatureQuantities} and the computation of $W^-$ from the appendix, we can compute the $L^2$ norms of $|\Ric|$ and $|\Riem|$.

To evaluate these integral norms on the 4-manifold parameterized by $(u,v,\theta_1,\theta_2)$, we use the parameterization
\begin{eqnarray}
	u,v\in[0,\infty) \quad \text{and} \quad \theta_1,\theta_2\in[0,\sqrt{8}\pi)
\end{eqnarray}
discussed in Section \ref{SubSecComputationAsymptotics}.
\begin{prop} \label{PropCurvNorms}
	The $L^2$ norms of the Ricci and Riemann tensors are
	\begin{eqnarray}
		\begin{aligned}
			&L^2(\Ric)\;=\;32\pi^2\frac{k^2}{1-k^2}, \quad
			L^2(W)\;=\;32\pi^2\frac{1+k^2}{1-k^2}, \\
			&\quad\quad\quad\quad
			L^2(\Riem)\;=\;32\pi^2\frac{1+3k^2}{1-k^2}.
		\end{aligned}
	\end{eqnarray}
\end{prop}
\begin{proof}
	In (\ref{EqnComputedRicSqVol}) we computed
	\begin{eqnarray}
		\begin{aligned}
			&|\Ric|^2dVol_4 \\
			&\quad\quad
			\;=\;\frac{32k^2uv}{\left(1+(1+k)u^2+(1-k)v^2\right)^3}\,du\wedge{d}v\wedge{}d\theta_1\wedge{}d\theta_2.
		\end{aligned}
	\end{eqnarray}
	Integrating along $\theta_1,\theta_2\in[0,\sqrt{8}\pi)$, we have
	\begin{eqnarray}
		\begin{aligned}
			&\quad\int_{N^4}|\Ric|^2dVol_4
			\;=\;8\pi^2\int_{\Sigma^2}\frac{32k^2uv}{\left(1+(1+k)u^2+(1-k)v^2\right)^3}
			\,du\wedge{d}v.
		\end{aligned}
	\end{eqnarray}
	Integrating $u$, $v$ from $0$ to $\infty$ gives
	\begin{eqnarray}
		\begin{aligned}
			&\quad\int_{N^4}|\Ric|^2dVol_4
			\;=\;\frac{32\pi^2k^2}{1-k^2}.
		\end{aligned}
	\end{eqnarray}
	Using the computation $|W^-|^2=96K_\Sigma{}^2$ of (\ref{EqnWeylTensorNorm}) and also using (\ref{EqnsAdoptedMetricSectional}) and (\ref{EqnMetricInUV}) to compute the volume form, we have
	\begin{eqnarray}
		\begin{aligned}
			&|W^-|^2
			\;=\;96M^2\left(\frac{-1+k\left((1+k)u^2-(1-k)v^2\right)}{\left(1+(1+k)u^2+(1-k)v^2\right)^3}\right)^2 \quad \text{and} \\
			&dVol_4=\frac{2uv}{M^2}\left(1+(1+k)u^2+(1-k)v^2\right)
			du\wedge{}dv\wedge{}d\theta_1\wedge{}d\theta_2
		\end{aligned}
	\end{eqnarray}
	Integrating $\theta_1$, $\theta_2$ along $[0,\sqrt{8}\pi)$ and $u,v$ along $[0,\infty)$ gives $L^2(|W^-|)=32\pi^2\frac{1+k^2}{1-k^2}$.
	
	The value of $\int|\Riem|^2dVol_4$ follows from the identity $|\Riem|^2=\frac{s^2}{6}+2|\cRic|^2+|W|^2$ and the fact that our manifolds are scalar-flat and half-conformally flat.
\end{proof}

The Chern-Gauss-Bonnet formula for the Euler class is $\chi(N^4)=\frac{1}{8\pi^2}\int\frac{R^2}{24}-\frac12|\cRic|^2+\frac14|W|^2$.
In some presentations the factor of $\frac14$ on the $|W|^2$ term is not present, which is due to norming $W$ as an operator $\bigwedge^2\rightarrow\bigwedge^2$ instead of as a tensor; see the discussion after Lemma \ref{LemmaRicEigenform} in the Appendix.
Using the $L^2$ norms of Proposition \ref{PropCurvNorms} we see immediately that $\chi(N^4)=1$, as expected.

We remark that the signature of these manifolds is zero, but $\int|W^+|^2-|W^-|^2\ne0$.
As a result, Lemma \ref{PropCurvNorms} can be used to compute $\eta$-invariants of various squashed 3-spheres.
We do not pursue this further however.

\section{Three kinds of Blowdown} \label{SectionBlowDownLimits}

The asymptotic geometry of open manifolds, including tangent cones at infinity and blowdown limits, are important in the study of open manifolds.
Our investigation of the generalized Taub-NUTs ends with an examination of their blowdown objects, where we find some surprises.

Given a metric $g$ on a complete $4$-manifold $N^4$, a Gromov-Hausdorff blowdown limit (colloquially known as {\it a tangent cone at infinity}) is a Gromov-Hausdorff limit of the manifold $N^4$ with metric $\epsilon_i^2g$ as $\epsilon_i\searrow0$.
In general such limits need not exist, and when they exist they need not be unique, and need not even be manifolds.

In this paper the objects have, for the most part, cubic volume growth and quadratic curvature decay $|\Riem|=O(r^{-2})$.
Therefore we expect limits to exist, and by computation we find they are unique.
The tangent cones at infinity are collapsed as expected, but they are not necessarily 3-dimensional as one might expect.

By (\ref{EqnMetricConvBlowDown}) we see the polytope metric itself converges uniquely under blowdown, but the situation on the full instanton is more complex, as there is a complicating geometric issue.
Level-sets of the distance function are spheres.
The collapsing field foliates the Hopf tori on these spheres, and this field might be rational or irrational.
Since $\mathcal{X}^1$ and $\mathcal{X}^2$ are the principle rotations on the level-sets and since the collapsing field is $\mathcal{X}=(1-k)\mathcal{X}^1-(1+k)\mathcal{X}^2$, this Hopf foliation is rational if and only if $k$ is rational.

In the rational case, the spherical level-sets converge down to $\mathbb{S}^2$ with up to two orbifold points, and in the irrational case the spherical level-sets converge down to a line segment (see Example 1.4 (continued) and Example 2.1 on pg 326 of \cite{CG1}).
Therefore when $k$ is rational the instanton blows down to a complete 3-dimensional stratified orbifold, and when $k$ is irrational it blows down to a 2-manifold with boundary. 

Because this behavior is rather pathological, we choose to modify the usual blowdown process in order to obtain better behavior in the limits.
The three kinds of blowdown are as follows.
The first is the usual Gromov-Hausdorff blowdown of $(N^4,g)$ itself, which we have just described.
The second and third kinds of blowdown, which we call ``generalized blowdowns,'' eliminate the pathologies arising from possibly irrational collapsing fields.
The second kind of blowdown is performed by ``unwrapping'' the torus fibers---this is just taking the interior of the polytope crossed with $\mathbb{R}^2$ instead of with the torus, and simply declining to apply the Delzant gluing process on the boundary.
So a new 4-manifold (not geodesically complete) exists with $\mathbb{R}^2$ fibers instead of torus fibers over each point of the polytope.
Then we take the blowdown limit of this object.
The metric has a zero eigen-direction on the $\mathbb{R}^2$ fibers.
We throw this direction away so now we have an $\mathbb{R}^1$-bundle over the polytope.
Finally we compactify the line fibers into circle fibers and so obtain a stratified limiting 3-conifold.
This conifold is an orbifold precisely when $k$ is rational.

For the third kind of ``blowdown'' limit, we take the blowdown on the polytope $\Sigma^2$ itself without regarding it as part of a larger 4-manifold.
This blowdown converges uniquely to a 2-manifold with boundary, and has a Riemannian metric with a curvature singularity at the origin.
Still, it has two well-defined momentum functions, so we can artificially perform the Delzant construction and still produce an honest 4-dimensional manifold with a Riemannian metric that has a point-like curvature singularity.

\subsection{Metric and coordinate convergence under blowdown}

From (\ref{EqnsAdoptedMetricSectional}) the polytope metric for the generalized Taub-NUTs is
\begin{eqnarray}
	g_\Sigma=\frac{2}{M}\left(1+(1+k)u^2+(1-k)v^2\right)\left(du^2+dv^2\right).
\end{eqnarray}
The 4-metric is $g_4=g_\Sigma+G^{ij}d\theta_i\otimes{d}\theta_j$ where $G^{ij}$ is the matrix from (\ref{EqnMetricInUV}).
Scaling the coordinates by setting $u=\sqrt[4]{M/2}\,\tilde{u}$, $v=\sqrt[4]{M/2}\,\tilde{v}$, the polytope metric becomes
\begin{eqnarray}
	\begin{aligned}
	g_\Sigma
	&=\left(\sqrt{\frac{2}{M}}+(1+k)\tilde{u}^2+(1-k)\tilde{v}^2\right)
	\left(d\tilde{u}^2+d\tilde{v}^2\right),
	\end{aligned}
\end{eqnarray}
and we can send $M\rightarrow\infty$.
Both metrics $g_\Sigma$ and $g_4$ converge.
We get
\begin{eqnarray}
	\begin{aligned}
		&g_\Sigma
		=\left((1+k)\tilde{u}^2+(1-k)\tilde{v}^2\right)\left(d\tilde{u}^2+d\tilde{v}^2\right), \\
		&G^{ij}\;=\;
		\frac{\frac12\tilde{u}^2\tilde{v}^2(\tilde{u}^2+\tilde{v}^2)}{(1+k)\tilde{u}^2+(1-k)\tilde{v}^2}
		\left(\begin{array}{cc}
		(1+k)^2  & 1-k^2 \\
		1-k^2 & (1-k)^2
		\end{array}\right). \label{EqnMetricConvBlowDown}
	\end{aligned}
\end{eqnarray}
Note that $\det(G^{ij})=0$.
Its zero eigenvector is $\vec{v}=(1-k)\frac{\partial}{\partial\theta_1}-(1+k)\frac{\partial}{\partial\theta_2}$ and its eigenvector of eigenvalue $\frac{\tilde{u}\tilde{v}(\tilde{u}^2+\tilde{v}^2)((1+k)^2+(1-k)^2)}{4((1+k)\tilde{u}^2+(1-k)\tilde{v}^2)}$ is $\vec{v}=(1+k)\frac{\partial}{\partial\theta_1}+(1-k)\frac{\partial}{\partial\theta_2}$.
Setting $\tilde\theta=\frac{1}{2(1+k^2)}\left((1+k)\theta_1+(1-k)\theta_2\right)$ gives a 3-dimensional metric of
\begin{eqnarray}
	\begin{aligned}
		\quad&\quad\quad{}g
		=\left((1+k)\tilde{u}^2+(1-k)\tilde{v}^2\right)\left(d\tilde{u}^2+d\tilde{v}^2\right)
		+\frac{\frac12\tilde{u}^2\tilde{v}^2(\tilde{u}^2+\tilde{v}^2)}{(1+k)\tilde{u}^2+(1-k)\tilde{v}^2}(d{\tilde\theta})^2.
		\label{EqnLimiting3Metric}
	\end{aligned}
\end{eqnarray}
The question of parameterization of $\tilde\theta$ can be determined as follows.
The field $\partial/\partial\tilde\theta$ produces a subgroup of the torus of slope $(1+k)/(1-k)$.
If $k=m/n$ is rational then this subgroup has slope $m+n/m-n$, and so, for $\tilde\theta$ to obtain consistent values in the limit, it must have parameterization $\tilde\theta\in[0,2\sqrt{2(1+k^2)}\pi/LCM(m+n,m-n))$.
Certainly if $k$ becomes irrational then the parametrization vanishes and hence the Gromov-Hausdorff limit collapses to a 2-dimensional object with no $\tilde\theta$ variable.

In the limit the metric (\ref{EqnLimiting3Metric}) no longer produces smooth points at the coordinate axes unless $k=0$.
Imitating the analysis at the beginning of Section \ref{SubSecComputationAsymptotics} at the coordinate axes, we can fix a positive value of $\tilde{u}$ and examine the 2-dimensional submanifold given by varying $\tilde{v}$, $\tilde\theta$, to obtain the metric
\begin{eqnarray}
	\begin{aligned}
		&\quad\quad{}g
		=
		(1+k)\tilde{u}^2
		\left(
		\left(1+O(\tilde{v}^2)\right)d\tilde{v}^2
		+\frac{1}{2(1+k)^2}\tilde{v}^2(1+O(\tilde{v}^2)(d{\tilde\theta})^2
		\right).
	\end{aligned}
\end{eqnarray}
and we can fix the value of $\tilde{u}$ and examine the 2-dimensional submanifold given by varying $\tilde{u}$, $\tilde\theta$, to obtain the metric
\begin{eqnarray}
	\begin{aligned}
		&\quad\quad{}g
		=
		(1-k)\tilde{v}^2
		\left(
		\left(1+O(\tilde{u}^2)\right)d\tilde{u}^2
		+\frac{1}{2(1-k)^2}\tilde{u}^2(1+O(\tilde{u}^2)(d{\tilde\theta})^2
		\right).
	\end{aligned}
\end{eqnarray}
When $k$ is rational we therefore observe a cone angle of $2\pi(1+k^2)/(1+k)LCM(m+n,m-n)$ along the $\tilde{v}$-axis and a cone angle of$2\pi(1+k^2)/2(1-k)LCM(m+n,m-n)$ along the $\tilde{u}$-axis,
These expressions are rational, so we observe a stratified orbifold in the limit.

When the collapsing field is irrational, the Gromov-Hausdorff limit does {\it not} produce such an object.
Following the discussion in Example 1.4 (continued) and Example 2.1 on pg 326 of \cite{CG1}, when the collapsing direction is irrational it collapses the spheres to line segments.

The central observation of the first ``generalized blowdown'' is that the metric (\ref{EqnLimiting3Metric}) makes sense on a 3-dimensional conifold, whether or not it is the result of a Gromov-Hausdorff blowdown.
We simply declare (\ref{EqnLimiting3Metric}) to be a metric on $(\tilde{u},\tilde{v},\tilde\theta)$ where $\tilde{u},\tilde{v}\in[0,\infty)$ and we give $\tilde\theta$ the range $[0,\sqrt{8}\pi)$ whether or not this is the range inherited from the Gromov-Hausdorff blowdown.
Of course one may give $\tilde\theta$ any other range---this will affect the cone angles but can never make both cone angles equal to $2\pi$ unless $k=0$, so can never produce a smooth manifold for $k\ne0$.
In the irrational case, both cone angles (along the two axes) cannot be simultaneously made rational, so for irrational $k$ we can produce a variety of stratified conifolds depending on the parameterization chosen for $\tilde\theta$, but we can never produce an orbifold.

We clearly still have a Killing field $\tilde{\mathcal{X}}={\partial}/{\partial\tilde\theta}$, and we may take a Riemannian quotient along to obtain the quarter-plane polytope again.
Its sectional curvature is
\begin{eqnarray}
	K_\Sigma\;=\;
	k\frac{
	(1+k)\tilde{u}^2-(1-k)\tilde{v}^2}{\left((1+k)\tilde{u}^2+(1-k)\tilde{v}^2\right)^3}.
	\label{EqnKBlowDown}
\end{eqnarray}
We see an irremovable curvature singularity at the origin $(\tilde{u},\tilde{v})=(0,0)$.
It is also not difficult to compute the Ricci curvature of the conifold.
It is diagonal in $(\tilde{u},\tilde{v},\tilde\theta)$ coordinates, and is given by
\begin{eqnarray}
	\Ric{}_3=
	\left(\begin{array}{ccc}
	\frac{-4k}{(1+k)\tilde{u}^2+(1-k)\tilde{v}^2} & 0 & 0 \\
	0 & \frac{4k}{(1+k)\tilde{u}^2+(1-k)\tilde{v}^2} & 0 \\
	0 & 0 & \frac{2k(1+k)^2\tilde{u}^2\tilde{v}^2(\tilde{u}^2-\tilde{v}^2)}{\left((1+k)\tilde{u}^2+(1-k)\tilde{v}^2\right)^3}
	\end{array}\right)
\end{eqnarray}
One notices immediately that scalar curvature is not zero.

\subsection{The second generalized blowdown}

The final type of ``blowdown'' comes out of the recognition that, in either of the two other blowdown processes, we always obtain a metric polytope, and that any such polytope does indeed encode all metric, complex structure, and curvature information for some scalar-flat 4-dimensional instanton, whether or not this 4-dimensional object has anything to do with any form of convergence of other metrics.

Using the expression (\ref{EqnsAdoptedMetricSectional}) for the moment functions and plugging in $\tilde{u}$, $\tilde{v}$, then sending $M\rightarrow\infty$, we have $\tilde\varphi^1=\frac{1+k}{2}\tilde{u}^2\tilde{v}^2$, $\tilde\varphi^2=\frac{1-k}{2}\tilde{u}^2\tilde{v}^2$, and we see the two rescaled moment functions are multiples of each other.
This gives a single moment function, which we set to $\tilde\varphi^1=\frac12\tilde{u}^2\tilde{v}^2$.

To obtain a second moment function, we perform a very natural renormalization.
Consider the function $\tilde\varphi^2=-(1-k)\varphi^1+(1+k)\varphi^2=-\frac{1-k}{M}v^2+\frac{1+k}{M}u^2$, which clearly coincides with the 0-eigenvector of the scaled metric.
To counteract the fact that the eigenvalue is approaching 0, we artificially scale $\tilde\varphi^2$ by $M$, and in the limit obtain $\tilde\varphi^2=-(1-k)\tilde{v}^2+(1+k)\tilde{u}^2$.
This gives us a second moment function.
We note that this renormalization process is directly analogous to the coordinate renormalization of Cheeger-Gromov in the proof of Theorem 2.1 of \cite{CG1}.
\begin{prop}[Third type of blowdown]
	Blowing down the polytope, the limit is the quarter-plane in $(\tilde{u},\tilde{v})$-coordinates.
	It has natural commuting momentum functions $\tilde\varphi^1=\frac12\tilde{u}^2\tilde{v}^2$, $\tilde{\varphi^2}=-\frac12(1+k)\tilde{u}^2+\frac12(1-k)\tilde{v}^2$.
	These are moment functions on the singular, toric, scalar-flat 4-conifold with a quarter-plane polytope that has metric
	\begin{eqnarray}
		\begin{aligned}
		&g_4\;=\;g_\Sigma\;+\;G^{ij}d\theta_i\otimes{d}\theta_j, \quad where \\
		&g_\Sigma\;=\;\left((1+k)\tilde{u}^2+(1-k)\tilde{v}^2\right)\left(d\tilde{u}\otimes{d}\tilde{u}+d\tilde{v}\otimes{d}\tilde{v}\right) \\
		&(G^{ij})\;=\;\left(\begin{array}{cc}
			\frac{\tilde{u}^2\tilde{v}^2(\tilde{u}^2+\tilde{v}^2)}{(1+k)\tilde{u}^2+(1-k)\tilde{v}^2} & \frac{-2k\tilde{u}^2\tilde{v}^2}{(1+k)\tilde{u}^2+(1-k)\tilde{v}^2} \\
			\frac{-2k\tilde{u}^2\tilde{v}^2}{(1+k)\tilde{u}^2+(1-k)\tilde{v}^2} & \frac{(1+k)^2\tilde{u}^2+(1-k)^2\tilde{v}^2}{(1+k)\tilde{u}^2+(1-k)\tilde{v}^2}
		\end{array}\right)
	\end{aligned}
\end{eqnarray}
\end{prop}
\begin{proof}
	The transitions from quadratic normal to volumetric normal coordinates are
	\begin{eqnarray}
		\begin{aligned}
		\tilde{u}=\sqrt[4]{M}\sqrt{\sqrt{x^2+y^2}+y}, \quad 
		\tilde{v}=\sqrt[4]{M}\sqrt{\sqrt{x^2+y^2}-y}.
		\end{aligned}
	\end{eqnarray}
	We also scale the volumetric coordinates, setting $x=\frac{1}{\sqrt{M}}\tilde{x}$, $y=\frac{1}{\sqrt{M}}\tilde{y}$, and obtain
	\be
		\begin{aligned}
			&\tilde{u}=\sqrt{\sqrt{\tilde{x}^2+\tilde{y}^2}+\tilde{y}}, \quad  \tilde{v}=\sqrt{\sqrt{\tilde{x}^2+\tilde{y}^2}-\tilde{y}} \\
			&\frac12\left(\tilde{u}^2+\tilde{v}^2\right)=\sqrt{\tilde{x}^2+\tilde{y}^2}, \quad \frac12\left(\tilde{u}^2-\tilde{v}^2\right)=\tilde{y}, \quad\tilde{u}\tilde{v}\;=\;\tilde{x}.
		\end{aligned}
	\ee
	Now we send $M\rightarrow\infty$.	
	The two limiting moment functions, in terms of $(\tilde{x},\tilde{y})$, are
	\be
		\begin{aligned}
		&\tilde\varphi^1=\frac12\tilde{u}^2\tilde{v}^2=\frac12\tilde{x}^2, \\
		&\tilde\varphi^2=-\frac12(1-k)\tilde{v}^2+\frac12(1+k)\tilde{u}^2=\tilde{y}\,+\,k\sqrt{\tilde{x}^2+\tilde{y}^2}.
		\end{aligned}
	\ee
	This is obviously a map from the right half-plane in the $(\tilde{x},\tilde{y})$ system to the half-plane in $\tilde\varphi^1$-$\tilde\varphi^2$ coordinates.
	We have transitions
	\be
		A=\left(\begin{array}{cc}
		\tilde{x} & 0\\
		\frac{k\tilde{x}}{\sqrt{\tilde{x}^2+\tilde{y}^2}} & 1+\frac{k\tilde{y}}{\sqrt{\tilde{x}^2+\tilde{y}^2}}
		\end{array}\right)
	\ee
	Using (\ref{EqnsMetricAndSectionalOnPolytope}) the polytope metric in $(\tilde{x},\tilde{y})$ and in $(\tilde{u},\tilde{v})$ coordinates is
	\begin{eqnarray}
		\begin{aligned}
			&g_\Sigma\;=\;\frac{k\tilde{y}+\sqrt{\tilde{x}^2+\tilde{y}^2}}{\sqrt{\tilde{x}^2+\tilde{y}^2}}\left(d\tilde{x}\otimes{d}\tilde{x}+d\tilde{y}\otimes{d}\tilde{y}\right) \\
			&g_\Sigma\;=\;\left((1+k)\tilde{u}^2+(1-k)\tilde{v}^2\right)\left(d\tilde{u}\otimes{d}\tilde{u}+d\tilde{v}\otimes{d}\tilde{v}\right).
		\end{aligned}
	\end{eqnarray}
	The corresponding 4-manifold metric is $g_4=g_\Sigma+G^{ij}d\theta_i\otimes{d}\theta_j$.
\end{proof}

\section{The exceptional Taub-NUT} \label{SectionExceptionalTaubNUT}

Unfortunately the ``almost distance function'' $\widetilde{R}$ of section \ref{SectionGeneralizedTN}, so crucial for determining manifold asymptotics, cannot be used when $k=\pm1$.
Here we imitate the analysis of Section \ref{SectionGeneralizedTN} in the exceptional case, and find a new almost distance function that is adapted to the exceptional case.

\subsection{Coordinates}

The exceptional Taub-NUT has moment functions, in terms of the volumetric normal coordinates, given by
\begin{eqnarray}
	\varphi^1=\frac{1}{\sqrt{2}}\left(-y+\sqrt{x^2+y^2}\right)\,+\,\frac{\alpha}{2}x^2, \;
	\varphi^2=\frac{1}{\sqrt{2}}\left(y+\sqrt{x^2+y^2}\right).
\end{eqnarray}
Simultaneous scaling in the $(x,y)$ and $(\varphi^1,\varphi^2)$ coordinates allows us to adjust $\alpha$, and we take $\alpha=2\sqrt{2}$, which is $M=1$.
Then (\ref{EqnsMetricAndSectionalOnPolytope}) gives the polytope metric
\begin{eqnarray}
	g_\Sigma\;=\;\frac{1+2y+2\sqrt{x^2+y^2}}{\sqrt{x^2+y^2}}\;\left(dx\otimes{d}x+dy\otimes{d}y\right).
\end{eqnarray}
The transitions to $(u,v)$ coordinates are
\begin{eqnarray}
	u=\sqrt{\sqrt{x^2+y^2}+y}, \quad v=\sqrt{\sqrt{x^2+y^2}-y}.
\end{eqnarray}
In these coordinates we may express the moment variables and polytope metric:
\begin{eqnarray}
	\begin{aligned}
		&\varphi^1\;=\;\frac{1}{\sqrt{2}}v^2(1+u^2), \quad
		\varphi^2\;=\;\frac{1}{\sqrt{2}}u^2,  \\
		&g_\Sigma\;=\;\left(1+u^2\right)\left(du\otimes{d}u+dv\otimes{d}v\right).
		\label{EqnExceptTNMetricUV}
	\end{aligned}
\end{eqnarray}
The matrix $G^{ij}=\left<\mathcal{X}^i,\mathcal{X}^j\right>$ is
\begin{eqnarray}
	(G^{ij})
	\;=\;\frac{1}{1+u^2}\left(\begin{array}{cc}
		v^2\left((1+u^2)^2+u^2v^2\right)\;\; & u^2v^2 \\
		\\
		u^2v^2 & u^2
	\end{array}\right) \label{EqnExcepMetricPerp}
\end{eqnarray}
so we have reconstructed the instanton metric: $g_4=g_\Sigma+G^{ij}d\theta_i\otimes{d}\theta_j$.
It is easy to compute the polytope sectional curvature in $(u,v)$ coordinates using (\ref{EqnsMetricAndSectionalOnPolytope}):
\begin{eqnarray}
	K_\Sigma\;=\;-\frac{1-u^2}{(1+u^2)^3}.
\end{eqnarray}
Notice that $K_\Sigma=-1$ along the positive $v$ axis.
Thus the instanton $(N^4,g_4)$ has $|\Riem|=O(1)$ along all geodesic rays that map to this ray in the polytope; these rays in $N^4$ constitute a rotationally symmetric 2-dimensional submanifold that is totally geodesic, as it is the zero-set of one of the Killing fields.
Using (\ref{EqnExceptionalTaubNUTDist}) below, where we compute geodesics and the distance function $R$, this also implies that $K_\Sigma=O(R^{-2})$ along all other geodesics based at the origin.

\subsection{Distance functions and geodesic normal coordinates}

To find distance functions, we imitate the separation trick of Section \ref{SectionGeneralizedTN}.
Setting $S(u,v)=f(u)+h(v)$ and finding
\begin{eqnarray}
	1=|\nabla{S}|^2\;=\;\frac{(f_u)^2+(h_v)^2}{1+u^2}.
\end{eqnarray}
Choosing a parameter $\eta\in[0,\pi/2]$ we have
\begin{eqnarray}
	\begin{aligned}
		&(f_u)^2+(h_v)^2\;=\;\left(\cos^2(\eta)+u^2\right)\,+\,\sin^2(\eta)
	\end{aligned}
\end{eqnarray}
so separating into $f_u=\sqrt{\cos^2(\eta)+u^2}$, $h_v=\sin(\eta)$ and integrating gives
\begin{eqnarray}
	\begin{aligned}
		&\quad\quad
		S_\eta(u,v)=\frac{\cos^2\eta}{2}\left(\frac{u}{\cos\eta}\sqrt{1+\frac{u^2}{\cos^2\eta}}+\log\left(\frac{u}{\cos\eta}+\sqrt{1+\frac{u^2}{\cos^2\eta}}\right)\right) \\
		& \quad\quad\quad\quad\quad\quad +v\sin(\eta)
	\end{aligned}
\end{eqnarray}
For each $\eta$, the characteristics of $S_\eta$ provides one geodesic from the origin, found by solving $\frac{d\gamma}{dt}=\nabla{}S_\eta$ with $\gamma(0)=(0,0)$.
This gives the system
\begin{eqnarray}
	\frac{du}{dt}=\frac{\sqrt{\cos^2\eta+u^2}}{1+u^2}, \quad\quad \frac{dv}{dt}=\frac{\sin\eta}{1+u^2} \label{EqnExceptionalParamGeod}
\end{eqnarray}
which is already partially separated, and can be evaluated in closed form.
But first, following the process of Section \ref{SectionGeneralizedTN}, we find the unparameterized geodesic equation.
Eliminating the $t$ parameter from equations (\ref{EqnExceptionalParamGeod}), we have
\begin{eqnarray}
	\begin{aligned}
		&\frac{dv}{du}\;=\;\frac{\sin\eta}{\sqrt{\cos^2\eta+u^2}},\;\; \; or\; \\
		&v\;=\;\sin(\eta)\,\log\left(\frac{u}{\cos\eta}+\sqrt{1+\frac{u^2}{\cos^2\eta}}\right).
	\end{aligned} \label{EqnExceptionalTaubNUTUnparametrizedGeod}
\end{eqnarray}
Integrating the first equation in (\ref{EqnExceptionalParamGeod}) gives
\begin{eqnarray}
	\begin{aligned}
		R
		&=\frac12u\sqrt{\cos^2\eta+u^2}
		+\frac{2-\cos^2\eta}{2}\log\left(\frac{u}{\cos\eta}+\sqrt{1+\frac{u^2}{\cos^2\eta}}\right) \\
		&=\frac12u\sqrt{\cos^2\eta+u^2}+v\frac{1+\sin^2\eta}{2\sin\eta}
	\end{aligned} \label{EqnExceptionalTaubNUTDist}
\end{eqnarray}
and so we have recovered the distance function $R=R(u,v)$, which is evaluated explicitly by solving for $\eta=\eta(u,v)$ from (\ref{EqnExceptionalTaubNUTUnparametrizedGeod}), then plugging into (\ref{EqnExceptionalTaubNUTDist}).

We define ``almost polar coordinates'' $(\tilde{R},\tilde\eta)$ for the exceptional Taub-NUT metric by
\begin{eqnarray}
	\begin{aligned}
		&\tilde{R}=\frac12u^2+v,\;\;
		\tilde\eta=\tan^{-1}(\sqrt{2v}/u),
		\quad\text{with inverses}\quad \\
		&u=\sqrt{2}\cos(\tilde\eta)\sqrt{\tilde{R}},\;\;
		v=\sin^2(\tilde\eta)\tilde{R}
	\end{aligned} \label{EqnExceptionalTNAlmostTransitions}
\end{eqnarray}
This ``almost distance function'' on the exceptional Taub-NUT is not as precise as the almost distance function we found on the generalized Taub-NUTs.
It approximates the distance function $R$ only to within a factor of $\sqrt{2}$.
\begin{prop}[Almost distance function] \label{PropExceptionalTaubNutAlmostDistance}
	Let $R=R(u,v)$ be the distance function and let $\tilde{R}=\frac12u^2+v$ be the almost distance function.
	Then for sufficiently large $\tilde{R}$
	\begin{eqnarray}
		1\;\le\;\frac{R}{\tilde{R}}\;\le\;\sqrt{2}.
	\end{eqnarray}
\end{prop}
\begin{proof}
	Put $u=\sqrt{2}\cos(\tilde\eta)\tilde{R}^{1/2}$, $v=\sin^2(\tilde\eta)\tilde{R}$.
	Using (\ref{EqnExceptionalTaubNUTUnparametrizedGeod}) and (\ref{EqnExceptionalTaubNUTDist}) and using (\ref{EqnExceptionalTNAlmostTransitions}) to transition to almost polar coordinates gives
	\begin{eqnarray}
		\begin{aligned}
			&R=\tilde{R}\left[\cos\tilde\eta\sqrt{\frac{\cos^2\eta}{2\tilde{R}}+\cos^2\tilde\eta}
			\;+\;\sin^2\widetilde{\eta}\,\frac{1+\sin^2\eta}{2\sin\eta}\right], \\
			&\frac{\sin^2\tilde\eta}{\sin\eta}=\frac{1}{\tilde{R}}\left[\log\sqrt{\tilde{R}}+\log\left(\frac{\sqrt{2}\cos\tilde\eta}{\cos\eta}+\sqrt{\frac{1}{\tilde{R}}+\frac{2\cos^2\tilde\eta}{\cos^2\eta}}\right)\right].
		\end{aligned} \label{EqnTwoTransExTN}
	\end{eqnarray}
	From $\frac{1}{\sin\eta}\ge1$, we have $\frac{1+\sin^2\tilde\eta}{2}\ge1$ and $\frac{\cos^2\eta}{2\tilde{R}}+\cos^2\tilde\eta\ge\cos^2\tilde\eta$.
	Then
	\begin{eqnarray}
		\begin{aligned}
			\frac{R}{\tilde{R}}
			&\;=\;\cos\tilde\eta\sqrt{\frac{\cos^2\eta}{2\tilde{R}}+\cos^2\tilde\eta}\;+\;\sin^2\tilde\eta\,\frac{1+\sin^2\eta}{2\sin\eta} \\
			&\;\ge\;\cos^2\tilde\eta\;+\;{\sin^2\tilde\eta}\;=\;1.
		\end{aligned}
	\end{eqnarray}
	The upper bound $R/\tilde{R}\le\sqrt{2}$ is slightly more involved.
	We perform the estimate in two parts: first if $\cos\eta\le\frac{1}{\sqrt{2}}$ then $\sin\eta\ge\frac{1}{\sqrt{2}}$ and we have simply
	\begin{eqnarray}
		\begin{aligned}
			\frac{R}{\tilde{R}}
			&\;=\;\cos\tilde\eta\sqrt{\frac{\cos^2\eta}{2\tilde{R}}+\cos^2\tilde\eta}
			\;+\;\sin^2\tilde\eta\left(\frac{1+\sin^2\eta}{2\sin\eta}\right) \\
			&\;\le\;(1+\epsilon(\tilde{R}))\cos^2\tilde\eta
			\;+\;\sin^2\tilde\eta\left(\frac{1+\sin^2\eta}{\sqrt{2}}\right)\;\le\;\sqrt{2}
		\end{aligned} \label{IneqFirstGenIneq}
	\end{eqnarray}
	Then if $\cos\eta\ge\frac{1}{\sqrt{2}}$ we have $\sin\eta<\frac{1}{\sqrt{2}}$ so estimating $\frac{\sin^2\tilde\eta}{\sin\eta}$ is tougher.
	But using the second equation in (\ref{EqnTwoTransExTN}) we can estimate
	\begin{eqnarray}
		\begin{aligned}
			\frac{\sin^2\tilde\eta}{\sin\eta}
			&=\tilde{R}^{-1}\left[\log\sqrt{\tilde{R}}+\log\left(\frac{\sqrt{2}\cos\tilde\eta}{\cos\eta}+\sqrt{\frac{1}{\tilde{R}}+\frac{2\cos^2\tilde\eta}{\cos^2\eta}}\right)\right] \\
			&\le\tilde{R}^{-1}\left[\log\sqrt{\tilde{R}}+\log\left(2\cos\tilde\eta+\sqrt{\frac{1}{\tilde{R}}+4\cos^2\tilde\eta}\right)\right] \\
			&\le\epsilon(\tilde{R})+(1+\epsilon(\tilde{R}))\log\left(4\cos\tilde\eta\right)
			\;\le\;\tilde{R}^{-1} \log5 \\
		\end{aligned}
	\end{eqnarray}
	so that
	\begin{eqnarray}
		\begin{aligned}
			\frac{R}{\tilde{R}}
			&\;=\;\cos\tilde\eta\sqrt{\frac{\cos^2\eta}{2\tilde{R}}+\cos^2\tilde\eta}\;+\;\frac{\sin^2\tilde\eta}{\sin\eta}\frac{1+\sin^2\eta}{2} \\
			&\;\le\;\cos^2\tilde\eta\sqrt{1+\frac{\cos^2\eta}{2\tilde{R}\cos^2\tilde\eta}}
			+\frac{1}{\tilde{R}}\frac{1+\sin^2\eta}{2} \\
			&\;=\;(1+\epsilon(\tilde{R}))\,\cos^2\tilde\eta
			\,+\,\epsilon(\tilde{R}).
		\end{aligned} \label{IneqSecondGenIneq}
	\end{eqnarray}
	the right-hand side of (\ref{IneqSecondGenIneq}) is certainly smaller than $\sqrt{2}$ for large $\tilde{R}$.
	This, with (\ref{IneqFirstGenIneq}), gives the result.
\end{proof}

An important question is whether the exceptional Taub-NUT metric is complete, which is unaddressed elsewhere in the literature.
Having computed the distance function to the origin, we can show that indeed the metric is complete.
\begin{prop}
	The exceptional Taub-NUT metric is complete.
\end{prop}
\begin{proof}
	We first show that there are no critical points of the distance function $R$, where $R:N^4\rightarrow\mathbb{R}$ is the distance to point $(u,v,\theta_1,\theta_2)=(0,0,0,0)\in{}N^4$.
	Recall how $R=R(u,v)$ is evaluated: one uses (\ref{EqnExceptionalTaubNUTUnparametrizedGeod}) to determine $\eta=\eta(u,v)$ and then one plugs $u$, $v$, $\eta$ in to (\ref{EqnExceptionalTaubNUTDist}).
	The distance to the origin is $\mathcal{X}^1$, $\mathcal{X}^2$ invariant, so $R(u,v,\theta_1,\theta_2)=R(u,v)$.
	
	We use a bit of textbook first order PDE theory to establish smoothness of the distance function $R$.
	If the PDE coefficients are smooth (which they are) and characteristics do not cross, then solutions are smooth.
	But the unparameterized equation (\ref{EqnExceptionalTaubNUTUnparametrizedGeod}) gives the characteristics, and there is one characteristic for each $\eta\in[0,\pi/2]$.
	To see that characteristics do not cross, we take a partial derivative of (\ref{EqnExceptionalTaubNUTUnparametrizedGeod}) with respect to the parameter $\eta$ and see that there are no stable points.
	Setting
	\begin{eqnarray}
		\mathcal{F}(u,v,\eta)=\frac{v}{\sin(\eta)}-\log\left(\frac{u}{\cos\eta}+\sqrt{1+\frac{u^2}{\cos^2\eta}} \right) \label{EqnOnlyNeedLogTerm}
	\end{eqnarray}
	so that for fixed $\eta$, then $\{\mathcal{F}=0\}$ defines a characteristic, a computation gives
	\begin{eqnarray}
		\begin{aligned}
			&\frac{\partial\mathcal{F}}{\partial\eta}
			\;=\;
			-\frac{\cos\eta}{\sin\eta}\frac{v}{\sin\eta}
			-\frac{\sin\eta}{\cos\eta}\frac{
				\left(\frac{u}{\cos\eta}\right)
				+\left(1+\frac{u^2}{\cos^2\eta}\right)^{-\frac12}\frac{u^2}{\cos^2\eta}
			}{\frac{u}{\cos\eta}+\sqrt{1+\frac{u^2}{\cos^2\eta}}}
		\end{aligned}
	\end{eqnarray}
	which is strictly negative in the range $\eta\in(0,\pi/2)$.
	Therefore characteristics do not cross in this range.
	To see that characteristics do not cross on the axes $\eta=0$, $\eta=\pi/2$, one looks at (\ref{EqnExceptionalTaubNUTUnparametrizedGeod}) directly
	\begin{eqnarray}
		v\;=\;\sin(\eta)\,\log\left(\frac{u}{\cos\eta}+\sqrt{1+\frac{u^2}{\cos^2\eta}}\right)
	\end{eqnarray}
	and sees that if $\eta\ne0,\pi/2$, one never has $u$ or $v$ equal to zero unless $(u,v)=(0,0)$.
	
	Because $R$ is smooth, there are not critical points of the distance function.
	This means that each metric ball $B_{R_0}=\{R<R_0\}$ is diffeomorphic to the Euclidean ball in the tangent space, via the exponential map.
	But $N^4$ is exhausted by the increasing union of such balls: $N^4=\bigcup_{R_0>0}B_R$.
	Therefore $N^4$ is complete.
\end{proof}

\subsection{Volume and curvature computations}

The ``almost ball'' of radius $S$, denoted $AB(S)$, is the set of points in $N^4$ with $\tilde{R}=v+\frac12u^2\le{}S$.
Using $\det{G}^{-1}=\frac14u^2v^2$, we obtain the 4-manifold volume form $dVol=\frac12uv(1+u^2)du\wedge{d}v\wedge{d}\theta^1\wedge{d}\theta^2$.
The almost-ball's volume is therefore
\begin{eqnarray}
	\begin{aligned}
		\quad\quad\;\; Vol\,AB(\tilde{R})
		&\;=\;8\pi^2\int_{0}^{\sqrt{2\tilde{R}}}
		\int_{0}^{\tilde{R}-\frac12u^2}\frac12uv(1+u^2)\,dv\,du \\
		&\;=\; \frac{\pi^2}{3}\left(\widetilde{R}^4+2\widetilde{R}^3\right).
	\end{aligned} \label{EqnExceptionalTaubNutAlmostBallVol}
\end{eqnarray}
\begin{prop}
	The exceptional Taub-NUT instanton has quartic asymptotic volume growth: $\Vol\,B(R)=O(R^4)$.
\end{prop}
\begin{proof}
	Combine (\ref{EqnExceptionalTaubNutAlmostBallVol}) with Proposition \ref{PropExceptionalTaubNutAlmostDistance}.
\end{proof}

Next we compute the Ricci potentials and the Ricci pseudo-volume form.
Using that $\sqrt{\mathcal{V}}=\sqrt{Det\,g^{-1}}=uv$ and using the formalism from Section \ref{SubSectionCurvatureQuantities} we have
\begin{eqnarray}
	\begin{aligned}
		\mathcal{R}^1
		&=\left<\nabla\log\mathcal{V},\,\nabla\varphi^1\right>
		=\sqrt{2}\frac{1+u^2+v^2}{1+u^2} \\
		\mathcal{R}^2
		&=\left<\nabla\log\mathcal{V},\,\nabla\varphi^2\right>
		=\sqrt{2}\frac{1}{1+u^2}.
	\end{aligned}
\end{eqnarray}
Taking exterior derivatives and using (\ref{EqnRicciVolForm}), we have
\begin{eqnarray}
	\begin{aligned}
		&|\Ric|^2dVol
		\;=\;\frac{16uv}{(1+u^2)^3}\,du\wedge{d}v\wedge{d}\theta^1\wedge{d}\theta^2, \\
		&|\Ric|^2
		\;=\;\frac{16}{(1+u^2)^4}.
	\end{aligned} \label{EqnExceptionalTaubNUTRic}
\end{eqnarray}
Integrating over the $(u,v)$ quarter-plane clearly gives an infinite value.

Notice that (\ref{EqnExceptionalTaubNUTRic}) gives that $|\Ric|^2=16$ along $u=0$.
Also notice, using (\ref{EqnExceptionalTaubNUTDist}), that along all other geodesics we have $|\Ric|=O(R^{-2})$.

\subsection{A scaled and an unscaled pointed limit}

\subsubsection{The blowdown}

Scaling the metric (\ref{EqnExceptTNMetricUV}) by $\frac{1}{M^4}$ and scaling $u$ and $v$ by ${M}$, we then send $M\rightarrow\infty$ to obtain the blowdown polytope metric
\begin{eqnarray}
	g_\Sigma\;=\;u^2\left(du^2+dv^2\right)
\end{eqnarray}
and the matrix $G^{ij}$ converges to
\begin{eqnarray}
	(G^{ij})
	=\frac12\left(\begin{array}{cc}
		v^2(u^2+v^2) & 0 \\
		\\
		0 & 0
	\end{array}\right).
\end{eqnarray}
After one throws away the $\theta_2$-direction, this gives a metric on a 3-manifold that is singular along $u=0$.

Finally we execute the ``third'' blowdown process.
During the blowdown process we scale the $\theta_2$ up by $M^2$.
In the limit we have a 4-manifold metric.
The matrix $G^{ij}=\left<\mathcal{X}^i,\mathcal{X}^j\right>$ becomes
\begin{eqnarray}
	(G^{ij})
	=\frac12\left(\begin{array}{cc}
		v^2(u^2+v^2) & v^2 \\
		\\
		v^2 & 1
	\end{array}\right).
\end{eqnarray}
The resulting instanton metric $g_4=g_\Sigma+G^{ij}d\theta_i\otimes{d}\theta_j$ is singular along the axis $u=0$.
This ``generalized'' blowdown has a 2-dimensional submanifold along which we have both topological and curvature singularities. 
The polytope sectional curvature is $K_\Sigma=-u^{-4}$, so the $v$-axis clearly holds singular curvature values.

\subsubsection{An unscaled pointed limit} \label{SubSubSecUnscaledLimit}

The exceptional Taub-NUT instanton $(N^4,g_4)$ has rays along which a sectional curvature equals $-1$.
A natural question is what happens when we take an {\it unscaled} pointed limit along such a ray.
We shall see that the resulting limit is the exceptional half-plane instanton.

The sectional curvature is constant along the $v$-axis, so we rechoose coordinates to center ourselves farther and farther along this axis.
For any $A>0$ set
\begin{eqnarray}
	\tilde{u}=u, \quad \tilde{v}=v-A.
\end{eqnarray}
The range of these coordinates is $\tilde{u}\in[0,\infty)$ and $\tilde{v}\in[-A,\infty)$ so in the limit the range is the entire half-plane.
In the Gromov-Hausdorff limit, the torus fibers actually become cylinders: the $\mathcal{X}^1$ direction becomes infinite, and since the field $\mathcal{X}^1$ itself becomes infinitely long, we must renormalize it.
For each $A$, choose new Killing fields
\begin{eqnarray}
	\widetilde{\mathcal{X}}^1=\frac{1}{2A}\left(\mathcal{X}^1-2\sqrt{2}A^2\mathcal{X}^2\right), \quad
	\widetilde{\mathcal{X}}^2=\sqrt{2}\mathcal{X}^2. \label{EqnsKillingFieldTransitions}
\end{eqnarray}
The polytope metric converges to
\begin{eqnarray}
	g_\Sigma=\left(1+\tilde{u}^2\right)\left(d\tilde{u}^2+d\tilde{v}^2\right)
\end{eqnarray}
and one can check directly that the matrix $G^{ij}$ converges to
\begin{eqnarray}
	(G^{ij})=\frac{1}{1+\tilde{u}^2}\left(\begin{array}{cc}
		(1+\tilde{u}^2)^2+4\tilde{u}^2\tilde{v}^2 & 2\tilde{u}^2\tilde{v} \\
		\\
		2\tilde{u}^2\tilde{v} & \tilde{u}^2
	\end{array}\right)
\end{eqnarray}
(we omit the tedious but straightforward computation).
Finally, notice that for the new Killing fields in (\ref{EqnsKillingFieldTransitions}) we have new moment functions $\tilde{\varphi}^1$, $\tilde{\varphi}^2$ defined up to a constant.
Choosing the constant appropriately, the transitions from old to new moment functions are $\tilde\varphi^1=\frac{1}{2A}\left(\varphi^1-A^2(1+2\sqrt{2}\varphi^2)\right)$, $\tilde{\varphi}^2=\sqrt{2}\varphi^2$.
These functions also converge, and in the limit as $A\rightarrow\infty$ we obtain
\be
\tilde{\varphi}^1\;=\;\tilde{v}+\tilde{v}\tilde{u}^2, \quad \tilde{\varphi}^2=\frac12\tilde{u}^2.
\ee
Comparing this data to the data laid out in Section \ref{SectionHalfPlaneInstanton} we see that this limiting Riemannian manifold is indeed the exceptional half-plane instanton (with the momentum variables switched).

\section{The exceptional half-plane instanton} \label{SectionHalfPlaneInstanton}

The exceptional half-plane instanton is given by
\begin{eqnarray}
	\varphi^1=\frac12x^2, \quad\quad \varphi^2=y+yx^2.
\end{eqnarray}
Using (\ref{EqnsMetricAndSectionalOnPolytope}) we obtain the polytope metric
\begin{eqnarray}
	g_\Sigma\;=\;\left(1+x^2\right)\left(dx\otimes{d}x+dy\otimes{d}y\right)
\end{eqnarray}
and the matrix $G^{ij}=\left<\mathcal{X}^i,\mathcal{X}^j\right>$ is
\begin{eqnarray}
	G^{ij}\;=\;
	\frac{1}{1+x^2}
	\left(\begin{array}{cc}
		x^2 & 2x^2y \\
		\\
		2x^2y & (1+x^2)^2+4x^2y^2
	\end{array}\right).
\end{eqnarray}
One notices the formal similarity with the exceptional Taub-NUT instanton.
There are two differences: the domains of the variables, and the size of the torus fibers.
Notice that $G^{ij}$ for the half-plane and the exceptional Taub-NUT are substantively different: for instance $G^{22}$ is never zero, reflecting the fact that the field $\mathcal{X}^2$ has no zeros

The formal similarity between this metric and the exceptional Taub-NUT metric allows us to use all of the polytope formalism, except that the domain is now the half-plane instead of the quarter-plane.

We have Ricci potentials $\mathcal{R}^1=\frac{1}{1+x^2}$ and $\mathcal{R}^2=\frac{2y}{1+x^2}$.
Then from (\ref{EqnRicciVolForm}) the norm-square of Ricci curvature is
\begin{eqnarray}
	\begin{aligned}
		&|\Ric|^2dVol\;=\;\frac{16x}{(1+x^2)^3}\,dy\wedge{d}x\wedge{d}\theta^1\wedge{d}\theta^2 \\
		&|\Ric|^2\;=\;\frac{16}{(1+x^2)^4}.
	\end{aligned}
\end{eqnarray}
As with the exceptional Taub-NUT we have $|\Ric|=16$ along geodesics within the 2-dimensional submanifold given by $x=0$ (which is a totally geodesic submanifold), and we have $|\Ric|=O(R^{-2})$ along all other geodesics.

\begin{prop}
	The exceptional half-plane instanton is geodesically complete.
\end{prop}
\begin{proof}
	The half-plane polytope metric $g=(1+x^2)(dx^2+dy^2)$, $x\ge0$, in volumetric normal coordinates is formally identical to the exceptional Taub-NUT metric $g=(1+u^2)(du^2+dv^2)$, $u,v\ge0$, in quadratic normal coordinates.
	
	Therefore the distance function $R(x,y)$ for the exceptional half-plane will be formally identical to the distance function for the exception Taub-NUT, except that $\eta$ may take the range $\eta\in[0,\pi]$.
	After accounting for this, the proof is identical.
\end{proof}

A final question of interest is whether or not the exceptional half-plane instanton {\it is} the exceptional Taub-NUT.
They have different polytopes, but conceivably these bear a relationship to each other like the half-plane and quarter-plane polytopes for flat $\mathbb{C}\times\mathbb{C}$.
That is, perhaps they are the same Riemannian manifold, but one has two rotational fields whereas the other has one rotational and one translational field.

\begin{prop}
	The exceptional Taub-NUT instanton is not the exceptional half-plane instanton.
\end{prop}
\begin{proof}
	We prove that the if $\mathcal{X}$ is any Killing field on the exceptional Taub-NUT, then necessarily $\mathcal{X}$ is a constant-coefficient combination of $\mathcal{X}_1$ and $\mathcal{X}_2$.
	Any such vector field has a zero, whereas the exceptional half-plane has a vector field without a zero---its $\mathcal{X}^2$ field.
	Thus, after proving that any Killing field on the exceptional Taub-NUT is a linear combination of $\mathcal{X}^1$, $\mathcal{X}^2$ and therefore has a zero, we have proven that the exceptional Taub-NUT cannot be the exceptional half-plane instanton.
	
	Let $\mathcal{X}=A_i\mathcal{X}^i+B_i\nabla\varphi^i$ be any vector field; we shall compute the Lie derivative $\mathcal{L}_{\mathcal{X}}g$ where $g=G_\Sigma+G^{ij}d\theta\otimes{}d\theta_j$ is the exceptional Taub-NUT metric; here $g_\Sigma=(1+u^2)(du^2+dv^2)$ and $G^{ij}$ is given by (\ref{EqnExcepMetricPerp}).
	By ({\it{i}}) of Lemma \ref{LemmaToricBasics} we have $\theta_i=G_{is}Jd\varphi^s$, and it is sometimes more convenient to express $g_4=g_\Sigma+G_{ij}Jd\varphi^i\otimes{}Jd\varphi^j$.
	
	We first consider the $B_i\nabla\varphi^i$ part of the vector field.
	Using $\mathcal{L}_X=i_xd+di_x$ we have
	\begin{eqnarray}
		\mathcal{L}_{B_i\nabla\varphi^i}d\varphi^j
		\;=\;d\left(B_iG^{ij}\right)
		\;=\;\frac{\partial(B_iG^{ij})}{\partial\varphi^s}d\varphi^s
	\end{eqnarray}
	and using ({\it{v}}) of Lemma \ref{LemmaToricBasics}, we have
	\begin{eqnarray}
		\begin{aligned}
			\mathcal{L}_{B_i\nabla\varphi^i}Jd\varphi^j
			&\;=\;i_{B_i\nabla\varphi^i}dJd\varphi^j
			\;=\;-i_{B_i\nabla\varphi^i}
			\left(\nabla\varphi^j(G_{st})d\varphi^s\wedge{}Jd\varphi^t\right) \\
			&\;=\;-B_i\nabla\varphi^j(G_{st})Jd\varphi^tG^{is}
			\;=\;B_iG_{st}\nabla\varphi^j(G^{is})Jd\varphi^t
		\end{aligned}
	\end{eqnarray}
	using the pseudo-K\"ahler condition, ({\it{ii}}) of Lemma \ref{LemmaToricBasics}, we have
	\begin{eqnarray*}
		\begin{aligned}
			&\mathcal{L}_{B_i\nabla\varphi^i}Jd\varphi^j
			\;=\;B_iG_{st}\nabla\varphi^j(G^{is})Jd\varphi^t
			\;=\;B_iG_{st}\nabla\varphi^s(G^{ij})Jd\varphi^t
			\;=\;B_i\frac{\partial{}G^{ij}}{\partial\varphi^t}Jd\varphi^t.
		\end{aligned}
	\end{eqnarray*}
	These two facts mean that the Lie derivative $\mathcal{L}_{B_i\nabla\varphi^i}$ preserves the vector space $span\{d\varphi^i\otimes{}d\varphi^j\}$ and also the vector space $span\{Jd\varphi^i\otimes{}Jd\varphi^j\}$.
	Therefore, in block form, we see
	\begin{eqnarray}
		\mathcal{L}_{B_i\nabla\varphi^i}g_4
		\;=\;\left(
			\begin{array}{c|c}
				C^{ij} &  \\
				\hline
				& D^{ij}
			\end{array}
		\right) \label{EqnBlockLie1}
	\end{eqnarray}
	where the off-diagonals are $2\times2$ zero matrices.

	Next we consider the Lie derivatives $\mathcal{L}_{A_iJ\nabla\varphi^i}d\varphi^j$ and $\mathcal{L}_{A_iJ\nabla\varphi^i}d\varphi^j$.
	Using $\mathcal{L}_{A_iJ\nabla\varphi^i}=i_{A_iJ\nabla\varphi^i}d+di_{A_iJ\nabla\varphi^i}$ and the fact that $(d\varphi^j)(J\nabla\varphi^i)=0$, we have
	\begin{eqnarray}
		\begin{aligned}
			\mathcal{L}_{A_iJ\nabla\varphi^i}d\varphi^j
			\;=\;0.
		\end{aligned}
	\end{eqnarray}
	Using the pseudo-K\"ahler condition we see
	\begin{eqnarray}
		\begin{aligned}
			\mathcal{L}_{A_iJ\nabla\varphi^i}Jd\varphi^j
			&\;=\;i_{A_iJ\nabla\varphi^i}dJd\varphi^j
			-d\left(A_iG^{ij}\right) \\
			&\;=\;-i_{A_iJ\nabla\varphi^i}
			\left(\nabla\varphi^j(G_{st})d\varphi^s\wedge{}Jd\varphi^t\right)
			-d\left(A_iG^{ij}\right) \\
			&\;=\;A_id(G^{ij})
			-d\left(A_iG^{ij}\right)
			\;=\;-G^{ij}dA_i
		\end{aligned}
	\end{eqnarray}
	Using the fact that the coefficient matrix $G^{ij}$ is invariant under $B_iJ\nabla\varphi^i$, we see that
	\begin{eqnarray}
		\begin{aligned}
			\mathcal{L}_{A_iJ\nabla\varphi^i}(g_4)
			&\;=\;G_{ij}\mathcal{L}_{A_iJ\nabla\varphi^i}\left(Jd\varphi^i\otimes{}Jd\varphi^j\right) \\
			&\;=\;-G_{ij}\left(G^{is}dA_s\otimes{}Jd\varphi^j+Jd\varphi^i\otimes\left(G^{js}dA_s \right) \right) \\
			&\;=\;
			-dA_j\otimes{}Jd\varphi^j
			-Jd\varphi^i\otimes{}dA_i.
		\end{aligned} \label{EqnLieSpecial}
	\end{eqnarray}
	Therefore, in block form, we see
	\begin{eqnarray}
		\mathcal{L}_{A_iJ\nabla\varphi^i}(g_4)
		\;=\;\left(\begin{array}{c|c}
			 & E^{ij} \\
			\hline
			E^{ji} & 
		\end{array}\right)
	\end{eqnarray}
	where the diagonals are blocks of $2\times2$ zero matrices.
	
	For $X=A_iJ\nabla\varphi^i+B_i\nabla\varphi^i$ to be Killing, we must have
	\begin{eqnarray}
		0\;=\;
		\mathcal{L}_{A_iJ\nabla\varphi^i+B_i\nabla\varphi^i}(g_4)
		\;=\;
		\mathcal{L}_{A_iJ\nabla\varphi^i}(g_4)
		+\mathcal{L}_{B_i\nabla\varphi^i}(g_4) \label{EqnBlockLie2}
	\end{eqnarray}
	From the block forms (\ref{EqnBlockLie1}) and (\ref{EqnBlockLie2}), both terms must individually be zero.
	
	Considering the second term, we show that $\mathcal{L}_{B_i\nabla\varphi^i}g_4=0$ means $B_i=0$.
	Note that the $g_\Sigma$ part of the metric $(1+u^2)(du^2+dv^2)$ has a single symmetry direction, $\partial/\partial{}v$, so the only candidate for $B_i\nabla\varphi^i$ is $B_i\nabla\varphi^i=\partial/\partial{}v$.
	To check that $\partial/\partial{}v$ does not fix $G^{ij}d\theta_i\otimes{}d\theta_j$, note that $\mathcal{L}_{\partial/\partial{}v}d\theta^i=d\mathcal{L}_{\partial/\partial{}v}\theta^i=0$, therefore $\mathcal{L}_{\partial/\partial{}v}\left(G^{ij}d\theta_i\otimes{}d\theta_j\right)=\frac{\partial{}G^{ij}}{\partial{v}}d\theta_i\otimes{}d\theta_j$.
	But the matrix $\frac{\partial{}G^{ij}}{\partial{v}}$ is not the zero matrix, as a glance at (\ref{EqnExcepMetricPerp}) can verify.
	Therefore $\mathcal{L}_{\partial/\partial{}v}\left(G^{ij}d\theta_i\otimes{}d\theta_j\right)\ne0$.
	We conclude that $\mathcal{L}_{B_i\nabla\varphi^i}g_4=0$ forces $B_i\nabla\varphi^i\ne\partial/\partial{v}$ and so $B_i\nabla\varphi^i=0$.
	
	Considering the second term, by (\ref{EqnLieSpecial}) if $\mathcal{L}_{A_iJ\nabla\varphi^i}(g_4)=0$ then we have
	\begin{eqnarray}
		0\;=\;dA_j\otimes{}Jd\varphi^j+Jd\varphi^i\otimes{}dA_i
	\end{eqnarray}
	which forces $dA_j=0$ meaning that $A_i=const$.
	We conclude that if $X=A_iJ\nabla\varphi^i+B_i\nabla\varphi^i$ is a Killing field, then it is a constant-coefficient combination of the Killing fields $J\nabla\varphi^1$, $J\nabla\varphi^2$, as claimed.
\end{proof}

\begin{appendix}
	
	\section{Appendix}
	
	We study the Weyl tensor on the generalized Taub-NUT metrics, and to a lesser extent, on toric K\"ahler 4-manifolds in general.
	Our analysis centers on a pair of two-forms: the form we call $\omega^-$ given by
	\begin{eqnarray}
	\omega^-\;=\;\frac{1}{\sqrt{det(G^{ij})}}\left(d\varphi^1\wedge{}d\varphi^2
	+Jd\varphi^1\wedge{}Jd\varphi^2\right) \label{EqnDefOfNegForm}
	\end{eqnarray}
	and the Ricci form $\rho=Ric(J\cdot,\cdot)$.
	We always have $\omega^-\in\bigwedge^-$ and $|\omega^-|=\sqrt{2}$.
	It is well-known that $\rho=\frac14s\omega+\rho_0$ where $\rho_0\in\bigwedge^-$, so in the scalar-flat case $\rho\in\bigwedge^-$.
	In this Appendix we prove a general result: $\omega^-$ is an eigenform for $W^-$ in the toric scalar flat case, and a specific result: on generalized Taub-NUTs the Ricci form $\rho\in\bigwedge^-$ is an eigenform of $W^-$.
	\begin{prop} \label{PropUniversalEigenform}
		Assume $(N^4,g_4,J,\mathcal{X}^1,\mathcal{X}^2)$ is a scalar-flat toric K\"ahler 4-manifold with moment functions $\varphi^1$ and $\varphi^2$.
		Then the length-$\sqrt{2}$ form $\omega^-$ of (\ref{EqnDefOfNegForm}) is an eigenform of the Weyl tensor:
		\begin{eqnarray}
		W^-(\omega^-)\;=\;2K_\Sigma\,\omega^-
		\end{eqnarray}
		where $K_\Sigma$ is the Gaussian curvature of the metric polygon $(\Sigma^2,g_\Sigma)$ associated to $N^4$.
	\end{prop}
	Below we show $\rho$ is orthogonal to $\omega^-$, meaning $\omega^-\wedge{}\rho=0$.
	Unless the toric manifold is Einstein, there is a 2-form $\rho^\perp$ of length $\sqrt{2}$ (unique up to sign) with $\rho^\perp\in\bigwedge^-$ and $\bigwedge^-=span_{\mathbb{R}}\{\omega^-,\rho_0,\rho^\perp\}$.
	\begin{prop}
		Assume $(M^4,g_4,J,\mathcal{X}^1,\mathcal{X}^2)$ is a generalized Taub-NUT.
		Then both $\rho,\rho^\perp\in\bigwedge{}^-$ are eigenforms of $W^-$:
		\begin{eqnarray}
		\begin{aligned}
		&W^-(\rho)\;=\;-4K_\Sigma\,\rho \\
		&W^-(\rho^\perp)\;=\;2K_\Sigma\,\rho^\perp.
		\end{aligned}
		\end{eqnarray}
		The Weyl tensor has 2 distinct eigenvalues: $2K_\Sigma$ (double) and $-4K_\Sigma$.
		We have
		\begin{eqnarray}
		W^-\;=\;K_\Sigma\left(\omega^-\otimes\omega^-
		-4\frac{\rho}{|\rho|}\otimes\frac{\rho}{|\rho|}
		+\rho^\perp\otimes\rho^\perp\right)
		\end{eqnarray}
		and $|W^-|^2=24K_\Sigma^2$.
	\end{prop}
	The Weyl tensor on most scalar-flat toric 4-manifolds has 3 distinct eigenvalues, as one would expect; only on the Taub-NUTs does this reduce to 2 eigenvalues.
	
	Before moving on we establish a bit of notation, following \cite{Derd}.
	Any 2-tensor $\gamma\in\bigotimes{}^2T^*N^4$ can be regarded as a map $\gamma:\bigwedge^1\rightarrow\bigwedge^1$, and we may compose two such maps $\gamma$, $\epsilon$ using the convention
	\begin{eqnarray}
	(\gamma\epsilon){}_{ij}\;=\;\gamma_{is}\,\epsilon{}^s{}_j. \label{EqnTwoFormProdConvention}
	\end{eqnarray}
	Any tensor of the form $F\in\bigwedge^2\otimes\bigwedge^2$ is a map $F:\bigwedge^2\rightarrow\bigwedge^2$; we use the convention
	\begin{eqnarray}
	F(\zeta)_{ij} \;=\; \frac12F{}_{ijkl}\zeta^{lk} \label{EqnRiemConvention}
	\end{eqnarray}
	(notice the reversal of indices on $\zeta$).
	We sometimes use $r$ for the symmetric 2-tensor $\Ric$, so we can switch seamlessly between our formulas and those of \cite{Derd}.

	\subsection{General toric K\"ahler 4-manifolds}
	
	Our first lemma establishes the computational features on toric K\"ahler 4-manifolds we shall require, and give a short proof of each assertion.
	\begin{lem}[Toric K\"ahler relations.]\label{LemmaToricBasics}
		Let $N^4$ be a toric K\"ahler manifold (not necessarily scalar-flat) with action potentials $\varphi^1$ and $\varphi^2$ and angle variables $\theta_1$, $\theta_2$ so that $\{d\varphi^1,d\theta^1,d\varphi^2,d\theta^2\}$ is an oriented (but not orthonormal) frame.
		
		{\it{i}}. We have $d\varphi^i=-G^{ij}Jd\theta_j$. \\
		{\it{ii}}. We have the ``pseudo-K\"ahler'' relations
		\begin{eqnarray}
			\nabla\varphi^i(G^{jk})\;=\;\nabla\varphi^j(G^{ik}),
			\quad\text{and}\quad
			\frac{\partial}{\partial\varphi^i}G_{jk}
			\;=\;\frac{\partial}{\partial\varphi^j}G_{ik}.
		\end{eqnarray}
		{\it{iii}}. The integral leaves of the distribution $\{\nabla\varphi^1,\nabla\varphi^2\}$ are totally geodesic, and therefore have identical intrinsic and extrinsic sectional curvatures, which is $K_\Sigma$. \\
		{\it{iv}}. We have
		$\nabla\frac{\partial}{\partial\varphi^k}
		=\Gamma_{ij}^kd\varphi^j\otimes\frac{\partial}{\partial\varphi^k}
		+\Gamma_{ij}^kJd\varphi^j\otimes{}J\frac{\partial}{\partial\varphi^k}$
		where $\Gamma_{ij}^k=\frac12G^{sk}G_{ij,s}.$ \\
		{\it{v}}. The Hessians $\nabla^2\varphi^k$ are $J$-invariant.
		We have
		\begin{eqnarray}
			\begin{aligned}
			&\nabla\nabla\varphi^k
			\;=\;
			-\Gamma_{ij}^k\,d\varphi^i\otimes\nabla\varphi^j
			-\Gamma_{ij}^k\,Jd\varphi^i\otimes{}J\nabla\varphi^j, \\
			&\nabla\left(J\nabla\varphi^k\right)
			\;=\;
			-\Gamma_{ij}^k\,d\varphi^i\otimes{}J\nabla\varphi^j
			+\Gamma_{ij}^k\,Jd\varphi^i\otimes{}\nabla\varphi^j.
			\end{aligned} \label{EqnHessComp}
		\end{eqnarray}
		{\it{vi}}. The exterior derivatives $dJd\varphi^k$ are $dJd\varphi^k=-2\Gamma_{ij}^k\,d\varphi^i\wedge{}Jd\varphi^j$. \\
		{\it{vii}}. The K\"ahler form is
		\begin{eqnarray}
			\omega\;=\;-d\varphi^i\wedge{}d\theta_i
			\;=\;-G_{ij}\,d\varphi^i\wedge{}Jd\varphi^j.
		\end{eqnarray}
		{\it{viii}}. The covariant derivative of $\omega^-$ has the form
		\begin{eqnarray}
			\begin{aligned}
			\nabla\omega^-
			\;\in\;span\{Jd\varphi^1,Jd\varphi^2\}\otimes\bigwedge{}^-
			\end{aligned} \label{EqnsSpanNablaPM}
		\end{eqnarray}
	\end{lem}
	\begin{proof}
		We move down the list of items, providing a short proof for each. \\
		\underline{Proof of {\it{i}}.}
		Using $\nabla\varphi^i=J\frac{\partial}{\partial\theta_i}$ we compute
		\begin{eqnarray}
		G_{ij}Jd\varphi^j\left(\frac{\partial}{\partial\theta_k}\right)
		\;=\;G_{ij}d\varphi^i\left(\nabla\varphi^k\right)
		\;=\;G_{ij}G^{jk}\;=\;\delta_i^k
		\end{eqnarray}
		so we conclude that $d\theta_i=G_{ij}Jd\varphi^j$, as claimed. \\	
		\underline{Proof of {\it{ii}}.}
		The pseudo-K\"ahler relations are equivalent to the toric relations $[\nabla\varphi^i,\nabla\varphi^j]=0$.
		For the first relation we have
		\begin{eqnarray}
		\begin{aligned}
		\nabla\varphi^i(G^{jk})
		&\;=\;
		\left<\nabla_{\nabla\varphi^i}\nabla\varphi^j,\,\nabla\varphi^k\right>
		+\left<\nabla\varphi^j,\,\nabla_{\nabla\varphi^i}\nabla\varphi^k\right> \\
		&\;=\;
		\left<\nabla_{\nabla\varphi^j}\nabla\varphi^i,\,\nabla\varphi^i\right>
		+\left<\nabla\varphi^i,\,\nabla_{\nabla\varphi^j}\nabla\varphi^k\right>
		\;=\;\nabla\varphi^j(G^{ik}).
		\end{aligned}
		\end{eqnarray}
		For the second relations we use $\frac{\partial}{\partial\varphi^i}=G_{iu}\nabla\varphi^u$ and $dG_{jk}=-G_{js}G_{kt}dG^{st}$.
		Then
		\begin{eqnarray*}
			\begin{aligned}
			\frac{\partial}{\partial\varphi^i}G_{jk}
			&=-G_{iu}G_{js}G_{kt}\nabla\varphi^u\left(G^{st}\right)
			=-G_{iu}G_{js}G_{kt}\nabla\varphi^s\left(G^{ut}\right)
			=\frac{\partial}{\partial\varphi^j}G_{ik}.
			\end{aligned}
		\end{eqnarray*}
		\underline{Proof of {\it{iii}}.}
		To see total geodesy of the integral leaves, we compute in two ways
		\begin{eqnarray}
			\begin{aligned}
			&\left<\nabla_{\nabla\varphi^i}\nabla\varphi^j,\,J\nabla\varphi^k\right>
			\;=\;\left<\nabla_{J\nabla\varphi^k}\nabla\varphi^j,
			\,\nabla\varphi^i\right>, \quad \text{and} \\
			&\left<\nabla_{\nabla\varphi^i}\nabla\varphi^j,\,J\nabla\varphi^k\right>
			\;=\;
			\left<\nabla_{\nabla\varphi^j}\nabla\varphi^i,\,J\nabla\varphi^k\right>
			\;=\;
			\left<\nabla_{J\nabla\varphi^k}\nabla\varphi^i,\,\nabla\varphi^j\right>
		\end{aligned}
		\end{eqnarray}
		Summing the two equations gives $2\left<\nabla_{\nabla\varphi^i}\nabla\varphi^j,\,J\nabla\varphi^k\right>=J\nabla\varphi^k(G^{ij})$ which is zero because $J\nabla\varphi^k$ is Killing.
		We conclude that the second fundamental form is zero. \\
		\underline{Proof of {\it{iv}}.}
		This follows easily from the textbook formula for $\Gamma^k_{ij}$, using $G_{ij,s}=G_{is,j}$. \\
		\underline{Proof of {\it{v}}}		
		To see $J$-invariance of the Hessians, with any fields $X$, $Y$ we compute
		\begin{eqnarray}
			\begin{array}{rll}
			\left<\nabla_{JX}\nabla\varphi^k,\,JY\right>
			&\;=\;-\left<\nabla_{JX}J\nabla\varphi^k,\,Y\right>
			\quad\quad
			& \text{Constancy of $J$} \\
			& \;=\;\left<\nabla_{Y}J\nabla\varphi^k,\,JX\right>
			& \text{$J\nabla\varphi^k$ is Killing} \\
			& \;=\;\left<\nabla_{Y}\nabla\varphi^k,\,X\right>
			& \text{Constancy of $J$} \\
			& \;=\;\left<\nabla_{X}\nabla\varphi^k,\,Y\right>
			& \text{Symmetry of $Hess(\varphi^k)$}
		\end{array}
		\end{eqnarray}
		The formula for $\nabla^2\varphi^i$ follows directly from $\left<\frac{\partial}{\partial\varphi^i},\nabla\varphi^k\right>=\delta^k_i$ and the fact that $\nabla\frac{\partial}{\partial\varphi^i}$. \\
		\underline{Proof {\it{vi}}}.
		Using the computation for $\nabla{}J\nabla\varphi^k$, we have
		\begin{eqnarray}
			\begin{aligned}
			dJd\varphi^k
			&\;=\;Alt(\nabla{}J\nabla\varphi^k) \\
			&\;=\;
			-\left(\Gamma^k_{[ij]}d\varphi^i\otimes{}Jd\varphi^j \right)
			+\left(\Gamma^k_{[ij]}Jd\varphi^i\otimes{}d\varphi^j \right) \\
			&\;=\;-2\Gamma^k_{ij}\,d\varphi^i\wedge{}Jd\varphi^j.
			\end{aligned}
		\end{eqnarray}
		\underline{Proof {\it{vii}}}.
		Surely $\omega=-G_{ij}d\varphi^i\wedge{}Jd\varphi^j\in\bigwedge^+$.
		Using $Jd\varphi^i=G^{is}d\theta_s$ we compute
		\begin{eqnarray}
		\begin{aligned}
		\omega(\cdot,J\cdot)
		&\;=\;G_{ij}d\varphi^i\otimes{}d\varphi^j+G_{ij}Jd\varphi^i\otimes{}Jd\varphi^j \\
		&\;=\;G_{ij}d\varphi^i\otimes{}d\varphi^j+G^{ij}d\theta_i\otimes{}d\theta_j
		\end{aligned}
		\end{eqnarray}
		which is precisely the metric $g_4$. \\
		\underline{Proof of {\it{viii}} }.
		We first show that $*d\varphi^1\wedge{}d\varphi^2=-Jd\varphi^1\wedge{}Jd\varphi^2$.
		Because the integral leaves are Lagrangian, we have $span\{d\varphi^1,d\varphi^2\}\perp{}span\{Jd\varphi^1,Jd\varphi^2\}$.
		Because $|d\varphi^1\wedge{}d\varphi^2|^2=|Jd\varphi^1\wedge{}Jd\varphi^2|^2$ we conclude that $*(d\varphi^1\wedge{}d\varphi^2)=\pm{}Jd\varphi^1\wedge{}Jd\varphi^2$.	
		To establish the sign, note that an oriented frame is $(d\varphi^1,d\theta_1,d\varphi^1,d\theta_2)$.
		Because $d\theta_i=-G_{ij}Jd\varphi^j$ and $det(-G_{ij})>0$, the frame $(d\varphi^1,Jd\varphi^1,d\varphi^2,Jd\varphi^2)$ is oriented.
		Therefore $*(d\varphi^1\wedge{}d\varphi^2)=-Jd\varphi^1\wedge{}Jd\varphi^2$, and we conlcude $d\varphi^1\wedge{}d\varphi^2\mp{}Jd\varphi^1\wedge{}Jd\varphi^2\in\bigwedge{}^\pm$.
		
		Thus $\omega^-\in\bigwedge{}^-$.
		To see $|\omega^-|=\sqrt{2}$, we compute $\omega^-\wedge\omega^-=-2det(G^{st})^{-1}d\varphi^1\wedge{}Jd\varphi^1\wedge{}d\varphi^2\wedge{}Jd\varphi^2=-2dVol$. \\
		\underline{Proof of {\it{viii}}}.
		Rather than a tedious computation of $\nabla\omega^-$, we take a shortcut.
		We have
		\begin{eqnarray}
			\begin{aligned}
			\bigwedge{}^2
			&\;=\;span\left\{
			d\varphi^1\wedge{}d\varphi^2,\;\;
			d\varphi^1\wedge{}Jd\varphi^1,\;\;
			d\varphi^1\wedge{}Jd\varphi^2,\right. \\
			&\quad\quad\quad\quad\quad
			\left.d\varphi^2\wedge{}Jd\varphi^1,\;\;
			d\varphi^2\wedge{}Jd\varphi^2,\;\;
			Jd\varphi^1\wedge{}Jd\varphi^2\right\}
			\end{aligned}
		\end{eqnarray}
		Then we note that, by total geodesy of the $d\varphi^1$-$d\varphi^2$ leaves,
		\begin{eqnarray}
			\nabla_{\nabla\varphi^i}(d\varphi^1\wedge{}d\varphi^2)
			\;\in\;span\{d\varphi^1\wedge{}d\varphi^2\}.
		\end{eqnarray}
		Using this and the covariant-constanct of $J$ we have
		\begin{eqnarray}
			\nabla_{\nabla\varphi^i}(Jd\varphi^1\wedge{}Jd\varphi^2)
			\;\in\;span\{Jd\varphi^1\wedge{}Jd\varphi^2\}.
		\end{eqnarray}
		It now follows that $\nabla_{\nabla\varphi^i}\omega^-\in{}span\{ \varphi^1\wedge\varphi^2,\;J\varphi^1\wedge{}J\varphi^2\}.$
		But the bundle $\bigwedge{}^-$ is covariant-constant, and therefore
		\begin{eqnarray}
			\begin{aligned}
			\nabla_{\nabla\varphi^i}\omega^-
			&\;\in\;span\{proj_{\bigwedge{}^-}Jd\varphi^1\wedge{}Jd\varphi^2\} \\
			&\;=\;span\{d\varphi^1\wedge{}d\varphi^2+Jd\varphi^1\wedge{}Jd\varphi^2\}
			\;=\;span\{\omega^-\}.
			\end{aligned}
		\end{eqnarray}
		But since $\omega^-$ has constant length, we have that $\left<\nabla_{\nabla\varphi^i}\omega^-,\omega^-\right>=\frac12\nabla\varphi^i|\omega^-|^2=0$.
		We conclude, as claimed, that $\nabla_{\nabla\varphi^i}\omega^-=0$ and so
		\begin{eqnarray}
			\begin{aligned}
			&\nabla\omega^-
			\;\in\;span\{J\nabla\varphi^1,J\nabla\varphi^2\}\otimes\bigwedge{}^-.
			\end{aligned}
		\end{eqnarray}
	\end{proof}
	
	\begin{lem}[Quaterionic Relations]
		Assume $(N^4,g_4,J,\mathcal{X}^1,\mathcal{X}^2)$ is a scalar-flat toric K\"ahler 4-manifold; in particular $\rho\in\bigwedge^-$.
		Then $\omega^-\wedge\rho=0$.
		Assuming $\rho\ne0$, then, referencing the product given in (\ref{EqnTwoFormProdConvention}), the 2-form
		\begin{eqnarray}
		\rho^\perp
		=\frac{1}{|\rho|}\omega^-\rho
		\end{eqnarray}
		has $\rho^\perp\in\bigwedge{}^-$ and $|\rho^\perp|=\sqrt{2}$, and we have the quaterionic relations
		\begin{eqnarray}
		\omega^-\frac{\rho}{|\rho|}=\rho^\perp, \quad
		\frac{\rho}{|\rho|}\rho^\perp=\omega^-, \quad
		\rho^\perp\omega^-=\frac{\rho}{|\rho|}.
		\end{eqnarray}
		and $\omega^-\omega^-=-2Id$, $\rho\rho=-|\rho|^2Id$, $\rho^\perp\rho^\perp=-2Id$.
	\end{lem}
	\begin{proof}
		With $\omega^-$ a multiple of $d\varphi^1\wedge{}d\varphi^2+Jd\varphi^1\wedge{}Jd\varphi^2$ and using $\rho=d\mathcal{R}^i\wedge{}d\theta_i$ from \S\ref{SubSectionCurvatureQuantities}, immediately $\omega^-\wedge\rho=0$.
		The quaterionic relations follow from the well-known fact that $span_{\mathbb{R}}\{Id\}\oplus\bigwedge^-$ is algebraically isomorphic to the quaternions.
	\end{proof}
	
	\begin{lem} \label{LemmaRicMap}
		In the scalar-flat case, the (symmetric) Ricci tensor is anti-invariant under $\omega^-$.
		Specifically
		\begin{eqnarray}
		r\omega^-\,+\,\omega^-r\;=\;0,
		\end{eqnarray}
		which is the same as $\Ric{}_{i}{}^s\omega_{sj}+\omega_{is}\Ric{}^{s}{}_{j}\;=\;0$.
	\end{lem}
	\begin{proof}
		We have $\Ric=\rho\omega$ where $\omega$ is the K\"ahler 2-form.
		Because $\omega\in\bigwedge^+$ and $\omega^-\in\bigwedge^-$ we certainly have that $\omega$ and $\omega^-$ commute: $\omega\omega^-=\omega^-\omega$.
		The quaterionic relations of the previous lemma give $\rho\omega^-=-\omega^-\rho$.
		Using these facts, we compute
		\begin{eqnarray}
			r\omega^-=\rho\omega\omega^-
			=\rho\omega^-\omega
			=-\omega^-\rho\omega\;=\;-\omega^-r.
		\end{eqnarray}
	\end{proof}
	\begin{lem}
		Referencing the product of (\ref{EqnRiemConvention}), in the scalar-flat case we have $\left(\Ric\KNP{}g\right)(\omega^-)=0$.
		As a consequence,
		\begin{eqnarray}
			\Riem(\omega^-)=W^-(\omega^-).
		\end{eqnarray}
	\end{lem}
	\begin{proof}
		We show that $(\Ric\KNP{}g)(\omega^-)=r\omega^-+\omega^-r$.
		To see this, we use $(\Ric\KNP{}g)_{ijkl}=r_{il}g_{jk}+r_{jk}g_{il}-r_{ik}g_{jl}-r_{jl}g_{ik}$ and compute
		\begin{eqnarray}
			\begin{aligned}
			(\Ric\KNPp{}g)(\omega^-)_{ij}
			&\;=\;
			\frac12\left(r_{il}g_{jk}+r_{jk}g_{il}
			-r_{ik}g_{jl}-r_{jl}g_{ik}\right)(\omega^-){}^{lk} \\
			&\;=\;\frac12\left(r_{il}(\omega^-){}^l{}_j+r_{jk}(\omega^-)_i{}^k
			-r_{ik}(\omega^-)_j{}^k-r_{jl}(\omega^-){}^l{}_i\right) \\
			&\;=\;r_{il}(\omega^-)^l{}_j+(\omega^-)_i{}^lr_{lj}
			\;=\;r\omega^-+\omega^-r\;=\;0.
			\end{aligned}
		\end{eqnarray}
		The rest follows from the Riemann tensor decomposition in the scalar-flat K\"ahler case: $\Riem=\frac12\left(\Ric\KNP{}g\right)+W^-$.
	\end{proof}
	\begin{lem} \label{LemmaAppBasics}
		Assume $(N^4,g_4,J,\mathcal{X}^1,\mathcal{X}^2)$ is a scalar-flat, toric K\"ahler 4-manifold.
		Then the form $\omega^-$ of (\ref{EqnDefOfNegForm}) is an eigenform of both $\Riem$ and $W^-$.
		We have
		\begin{eqnarray}
		\begin{aligned}
		&\Riem(\omega^-)\;=\;
		W^-(\omega^-)
		\;=\;2K_\Sigma\omega^-.
		\end{aligned}
		\end{eqnarray}
	\end{lem}
	\begin{proof}		
		The $J$-invariance of $\Riem$ means $\Riem(\nabla\varphi^1,\nabla\varphi^2,\cdot,\cdot)=\Riem(J\nabla\varphi^1,J\nabla\varphi^2,\cdot,\cdot)$, and so $\Riem(\omega^-)=2det(G^{st})^{-1/2}\Riem(\nabla\varphi^1,\nabla\varphi^2,\cdot,\cdot)$.
		By the total geodesy of the distribution $span\{\nabla\varphi^1,\nabla\varphi^2\}$ we have $\Riem(\nabla\varphi^1,\nabla\varphi^2)\nabla\varphi^i\in{}span\{\nabla\varphi^1,\nabla\varphi^2\}$.
		Because $J$ is covariant-constant, we have $\Riem(\nabla\varphi^1,\nabla\varphi^2)\mathcal{X}_i\in{}span\{\mathcal{X}_1,\mathcal{X}_2\}$.
		Thus $\Riem(\nabla\varphi^1,\nabla\varphi^2,\nabla\varphi^i,\mathcal{X}_j)=0$ for any $i,j$.
		
		Therefore the only non-zero terms in $\Riem(\omega^-)$ are multiples of $\Riem(\nabla\varphi^1,\nabla\varphi^2,\nabla\varphi^2,\nabla\varphi^1)$ and $\Riem(\nabla\varphi^1,\nabla\varphi^2,\mathcal{X}_2,\mathcal{X}_1)$.
		By $J$-invariance again, we see $\Riem(\nabla\varphi^1,\nabla\varphi^2,\nabla\varphi^2,\nabla\varphi^1)=\Riem(\nabla\varphi^1,\nabla\varphi^2,\mathcal{X}_2,\mathcal{X}_1)=K_\Sigma\cdot\det(G^{st})$.
		Combining terms in the tensor $\frac12\Riem_{ij}{}^{kl}\omega^-{}_{lk}$ we therefore obtain
		\begin{eqnarray}
			&&\Riem(\omega^-)\;=\;2K_\Sigma\,\omega^-.
		\end{eqnarray}
		In the scalar-flat toric case the Riemann tensor decomposes as $\Riem=\frac12(\Ric\KNP{}g)+W^-$.
		We have shown above that $(\Ric\KNP{}g)(\omega^-)=0$, so $W^-(\omega^-)=2K_\Sigma\omega^-$ as claimed.
	\end{proof}

	\subsection{Specialization to the Taub-NUT metrics}
	
	To explore the Weyl tensor further, we the Derdzinski's framework of \cite{Derd}.
	In the K\"ahler case Derdzhinski has told us $W^+=\frac{s}{24}\left(3\omega\otimes\omega-2Id_{\bigwedge^+}\right)$ where $\omega$ is the K\"ahler form.
	Using (29) of \cite{Derd} for $W^-$ we have
	\begin{eqnarray}
		W^-\;=\;\frac12\left(
		\lambda^-\,\omega^-\otimes\omega^-
		\,+\,\mu^-\,\eta^-\otimes\eta^-
		\,+\,\nu^-\,\theta^-\otimes\theta^-
		\right)
	\end{eqnarray}
	where $\omega^-,\eta^-,\theta^-\in\bigwedge^-$ are the eigenforms of $W^-$ of length $\sqrt{2}$ and $\lambda^-,\mu^-,\nu^-\in\mathbb{R}$ are the corresponding eigenvalues.
	The forms $\omega^-$, $\eta^-$, $\theta^-$ are the length-$\sqrt{2}$ eigenforms of $W^-$ with eigenvales $\lambda$, $\mu$, $\nu$.
	We have the quaternionic relations $\omega^-\eta^-=\theta^-$ and cyclic permutations.
	From Proposition \ref{LemmaAppBasics} we certainly have $\lambda^-=2K_\Sigma.$
	The governing equations are (32) of \cite{Derd}
	\begin{eqnarray}
	\begin{array}{rlll}
	\nabla\omega^-&\;=\; &\;\;\;\;\,b\otimes\eta^-&\,-\;\;c\otimes\theta^- \\
	\nabla\eta^-&\;=\;-b\otimes\omega^-& &\,+\;\;a\otimes\theta^- \\
	\nabla\theta^-&\;=\;\;\;\;c\otimes\omega^-&\,-\;a\otimes\eta^-
	\end{array} \label{EqnsDerivEForms}
	\end{eqnarray}
	and (33) of \cite{Derd}
	\begin{eqnarray}
	\begin{aligned}
	&da\,+\,b\wedge{}c
	\;=\;\left(\lambda^--s/6\right)\omega^-\,+\,\left(\omega^-{}r+r\omega^-\right)/2 \\
	&db\,+\,c\wedge{}a
	\;=\;\left(\mu^--s/6\right)\eta^-\,+\,\left(\eta^-{}r+r\eta^-\right)/2 \\
	&dc\,+\,a\wedge{}b
	\;=\;\left(\nu^--s/6\right)\theta^-\,+\,\left(\theta^-{}r+r\theta^-\right)/2
	\end{aligned} \label{EqnsCurvOfWedgePl}
	\end{eqnarray}
	where $s$ is scalar curvature.
	From ({\it{viii}}) of Proposition \ref{LemmaToricBasics} we certainly have
	\begin{eqnarray}
	b,\,c\;\in\;span{}_{\mathcal{R}}\left\{Jd\varphi^1,\,Jd\varphi^2\right\}.
	\end{eqnarray}
	
	At this point we are forced to recess from general considerations and conduct computations.
	Because of our specialization from toric metrics generally to the Taub-NUTs, there is no other way to proceed.
	In the generalized Taub-NUT case the 1-forms $a$, $b$, and $c$ are
	\begin{eqnarray}
		\begin{aligned}
		&a=F^{-1}\left[
		2My\cdot{}dx
		-2M\left(x+k\sqrt{x^2+y^2}\right)
		\cdot{}dy
		\right] \\
		&b=F^{-1}\left[(-2kMy\cdot{}Jdx
		+2M(kx+\sqrt{x^2+y^2})\cdot{}Jdy\right] \\
		&c=F^{-1}\left[Jdx-
		\frac{1}{y}\left(x+2M\sqrt{x^2+y^2}\left(x+k\sqrt{x^2+y^2}\right)\right)
		\cdot{}Jdy\right]
	\end{aligned} \label{EqnABCDerd}
	\end{eqnarray}
	where $F=\sqrt{x^2+y^2}\left(1+2M\left(kx+\sqrt{x^2+y^2}\right)\right)$.
	The length-$\sqrt{2}$ eigenforms are
	\begin{eqnarray}
		\begin{aligned}
		\omega^-, \quad\;\;
		\eta^-=\sqrt{2}\frac{\rho}{|\rho|}, \quad\;\;
		\theta^-=\rho^\perp.
		\end{aligned}
	\end{eqnarray}
	\begin{lem}[The Ricci form] \label{LemDerivRic}
		Assume $(N^4,g_4,J,\mathcal{X}^1,\mathcal{X}^2)$ is a generalized Taub-NUT.
		Letting $\rho$ be the Ricci form, we have covariant derivative
		\begin{eqnarray}
		\nabla\left(\sqrt{2}\frac{\rho}{|\rho|}\right)\;=\;-b\otimes\omega^-
		\,+\,a\otimes\rho^\perp \label{EqnDerivRho}
		\end{eqnarray}
		where $a$ and $b$ are given by (\ref{EqnABCDerd}).
		The Laplacian of $\rho$ is $\triangle\rho=8K_\Sigma\rho$.
	\end{lem}
	\begin{proof}
		Equation (\ref{EqnDerivRho}) follows from the formula (\ref{EqnSpecificRicciPots}) for the Ricci potentials, from which $\rho=\Ric(J\cdot,\cdot)$ can be found, along with the formula for $\Gamma^i_{jk}$.
		To compute $\triangle\rho$,
		\begin{eqnarray*}
			\begin{aligned}
				&(\triangle\rho)_{ij}\;=\;g^{kl}\rho_{ij,kl}
				\;=\;g^{kl}\frac{\partial}{\partial{x}^k}\frac{\partial}{\partial{x}^l}\rho_{ij} \\
				&\;\;
				-g^{kl}\frac{\partial}{\partial{}x^k}\left(\Gamma_{il}^s\rho_{sj}\right)
				-g^{kl}\frac{\partial}{\partial{}x^k}\left(\Gamma_{jl}^s\rho_{is}\right)
				-g^{kl}\Gamma_{il}^s\frac{\partial}{\partial{}x^k}\rho_{sj}
				-g^{kl}\Gamma_{jl}^s\frac{\partial}{\partial{}x^k}\rho_{is}
				-g^{kl}\Gamma^s_{kl}\frac{\partial}{\partial{}x^k}\rho_{ij} \\
				&\;\;
				+g^{kl}\Gamma^t_{ik}\Gamma^s_{tl}\rho_{sj}
				+g^{kl}\Gamma^t_{jk}\Gamma^s_{tl}\rho_{is}
				+2g^{kl}\Gamma^t_{ik}\Gamma^s_{jl}\rho_{st}
				+g^{kl}\Gamma^t_{kl}\Gamma^s_{it}\rho_{sj}
				+g^{kl}\Gamma^t_{kl}\Gamma^s_{jt}\rho_{is}.
			\end{aligned}
		\end{eqnarray*}
		Fully worked out, this expression has no fewer than 14976 terms with 384 derivative operations, so computer assistance is essential.
		Using (\ref{EqnSpecificRicciPots}), (\ref{EqnsAdoptedMetricSectional}), and (\ref{EqnMetricInUV}), a short Mathematica code provides the result.
	\end{proof}
	
	\begin{lem} \label{LemmaRicEigenform}
		For the generalized Taub-NUT metrics, we have $W^-(\rho)=-4K_\Sigma\rho$.
		As a consequence we have
		\begin{eqnarray}
		W^-\;=\;K_\Sigma\left(\omega^-\otimes\omega^-
		-4\frac{\rho}{|\rho|}\otimes\frac{\rho}{|\rho|}
		+\rho^\perp\otimes\rho^\perp \right)
		\end{eqnarray}
		and $|W^-|^2\;=\;24K_\Sigma^2$.
	\end{lem}
	\begin{proof}
		The fact that $W^-(\rho)=-4K_\Sigma\rho$ follows directly from the Bochner formula of \cite{Le1}, which is $\triangle\rho=-2W^-(\rho^-)+\frac13s\rho$.
		Because scalar curvature is zero, this gives $8K_\Sigma\rho=-2W^{-}(\rho)$.		
		The expression for $W^-$ now follows from the fact that $\omega^-$ and $\rho$ are eigenforms so $\rho^\perp$ must be the final eigenform.
		The fact that $W^-$ is trace free forces $W^-(\rho^\perp)=2K_\Sigma\rho^\perp$.
		
		The expression for $|W^-|^2$ follows from the fact that the three terms in parentheses are mutually orthogonal, combined with $|\omega^-\otimes\omega^-|=4$, $|4\rho\otimes\rho|^2=16|\rho|^2$, and $|\rho^\perp\otimes\rho^\perp|^2=4$.
	\end{proof}
	Because $\omega^-$, $\sqrt{2}\rho/|\rho|$, $\rho^\perp$ are orthogonal anti self-dual 2-forms of length $\sqrt{2}$, we have $Id_{\bigwedge^-}=\frac12\left(\omega^-\otimes\omega^-+2|\rho|^{-2}\rho\otimes\rho+\rho^\perp\otimes\rho^\perp\right)$ and so we have expression (\ref{EqnsWeylTensorComplete}):
	\begin{eqnarray}
		\begin{aligned}
			W^-
			&\;=\;K_\Sigma\left(-4|\rho|^{-2}\rho\otimes\rho+\omega\otimes\omega+\rho^\perp\otimes\rho^\perp\right) \\
			&\;=\;
			K_\Sigma\left(-6|\rho|^{-2}\rho\otimes\rho+2Id_{\bigwedge^-}\right).
		\end{aligned}
	\end{eqnarray}
	
	{\bf Remark} In Theorem (\ref{LemmaRicEigenform}), following Derdzinski, we used the operator norm for $|W^-|^2$, where the operator $W^-:\bigwedge^-\rightarrow\bigwedge^-$ is described by (\ref{EqnRiemConvention}).
	This is not the standard tensor norm, but $\frac14$ times the standard tensor norm.
	Using the tensor norm
	\begin{eqnarray}
		|W^-|^2_{tensor}=W^-_{ijkl}W^-_{stuv}g^{is}g^{jt}g^{ku}g^{lv}, \label{EqnWeylTensorNorm}
	\end{eqnarray}
	we have that $|W^-|^2_{tensor}=96K_\Sigma^2$.
	This is an important point in Section \ref{SubSecEnergyComp}.

\end{appendix}

\end{document}